# Harold Jeffreys's Theory of Probability Revisited

**Christian P. Robert, Nicolas Chopin and Judith Rousseau**

*Abstract.* Published exactly seventy years ago, Jeffreys's *Theory of Probability* (1939) has had a unique impact on the Bayesian community and is now considered to be one of the main classics in Bayesian Statistics as well as the initiator of the objective Bayes school. In particular, its advances on the derivation of noninformative priors as well as on the scaling of Bayes factors have had a lasting impact on the field. However, the book reflects the characteristics of the time, especially in terms of mathematical rigor. In this paper we point out the fundamental aspects of this reference work, especially the thorough coverage of testing problems and the construction of both estimation and testing noninformative priors based on functional divergences. Our major aim here is to help modern readers in navigating in this difficult text and in concentrating on passages that are still relevant today.

## 1. INTRODUCTION

*The theory of probability makes it possible to respect the great men on whose shoulders we stand.*
H. Jeffreys, *Theory of Probability*, Section 1.6.

Few Bayesian books other than *Theory of Probability* are so often cited as a foundational text.[1] This book is rightly considered as the principal reference in modern Bayesian statistics. Among other innovations, *Theory of Probability* states the general principle for deriving noninformative priors from the sampling distribution, using Fisher information. It also proposes a clear processing of Bayesian testing, including the dimension-free scaling of Bayes factors. This comprehensive treatment of Bayesian inference from an objective Bayes perspective is a major innovation for the time, and it has certainly contributed to the advance of a field that was then submitted to severe criticisms by R. A. Fisher (Aldrich, 2008) and others, and was in danger of becoming a feature of the past. As pointed out by Zellner (1980) in his introduction to a volume of essays in honor of Harold Jeffreys, a fundamental strength of *Theory of Probability* is its affirmation of a unitarian principle in the statistical processing of all fields of science.

*Christian P. Robert is Professor of Statistics, Applied Mathematics Department, Université Paris Dauphine and Head of the Statistics Laboratory, Center for Research in Economics and Statistics (CREST), National Institute for Statistics and Economic Studies (INSEE), Paris, France e-mail: xian@ceremade.dauphine.fr. He was the President of ISBA (International Society for Bayesian Analysis) for 2008. Nicolas Chopin is Professor of Statistics, ENSAE (National School for Statistics and Economic Administration), and Member of the Statistics Laboratory, Center for Research in Economics and Statistics (CREST), National Institute for Statistics and Economic Studies (INSEE), Paris, France e-mail: nicolas.chopin@ensae.fr. Judith Rousseau is Professor of Statistics, Applied Mathematics Department, Université Paris Dauphine, and Member of the Statistics Laboratory, Center for Research in Economics and Statistics (CREST), National Institute for Statistics and Economic Studies (INSEE), Paris, France e-mail: rousseau@ensae.fr.*

[1]Among the "Bayesian classics," only Savage (1954), DeGroot (1970) and Berger (1985) seem to get more citations than Jeffreys (1939, 1948, 1961), the more recent book by Bernardo and Smith (1994) coming fairly close. The homonymous *Theory of Probability* by de Finetti (1974, 1975) gets quoted a third as much (*Source: Google Scholar*).

For a 21st century reader, Jeffreys's *Theory of Probability* is nonetheless puzzling for its lack of formalism, including its difficulties in handling improper priors, its reliance on intuition, its long debate about the nature of probability, and its repeated attempts at philosophical justifications. The title itself is misleading in that there is absolutely no exposition of the mathematical bases of probability theory in the sense of Billingsley (1986) or Feller (1970): "Theory of Inverse Probability" would have been more accurate. In other words, the style of the book appears to be both verbose and often vague in its mathematical foundations for a modern reader.[2] (Good, 1980, also acknowledges that many passages of the book are "obscure.") It is thus difficult to extract from this dense text the principles that made *Theory of Probability* the reference it is nowadays. In this paper we endeavor to revisit the book from a Bayesian perspective, in order to separate foundational principles from less relevant parts.

This review is neither a historical nor a critical exercise: while conscious that *Theory of Probability* reflects the idiosyncrasies both of the scientific achievements of the 1930's—with, in particular, the emerging formalization of Probability as a branch of Mathematics against the ongoing debate on the nature of probabilities—and of Jeffreys's background—as a geophysicist—, we aim rather at providing the modern reader with a reading guide, focusing on the pioneering advances made by this book. Parts that correspond to the lack (at the time) of analytical (like matrix algebra) or numerical (like simulation) tools and their substitution by approximation devices (that are not used any longer, even though they may be surprisingly accurate), and parts that are linked with Bayesian perspectives will be covered fleetingly. Thus, when pointing out notions that may seem outdated or even mathematically unsound by modern standards, our only aim is to help the modern reader stroll past them, and we apologize in advance if, despite our intent, our tone seems overly presumptuous: it is rather a reflection of our ignorance of the current conditions at the time since (to borrow from the above quote which may sound itself somehow presumptuous) we stand respectfully at the feet of this giant of Bayesian Statistics.

The plan of the paper follows *Theory of Probability* linearly by allocating a section to each chapter of the book (Appendices are only mentioned throughout the paper). Section 10 contains a brief conclusion. Note that, in the following, words, sentences or passages quoted from *Theory of Probability* are written in italics with no precise indication of their location, in order to keep the style as light as possible. We also stress that our review is based on the third edition of *Theory of Probability* (Jeffreys, 1961), since this is both the most matured and the most available version (through the last reprint by Oxford University Press in 1998). Contemporary reviews of *Theory of Probability* are found in Good (1962) and Lindley (1962).

## 2. CHAPTER I: FUNDAMENTAL NOTIONS

*The posterior probabilities of the hypotheses are proportional to the products of the prior probabilities and the likelihoods.*
H. Jeffreys, *Theory of Probability*, Section 1.2.

The first chapter of *Theory of Probability* sets general goals for a *coherent* theory of induction. More importantly, it proposes an axiomatic (if slightly tautological) derivation of prior distributions, while justifying this approach as coherent, compatible with the *ordinary process of learning* and allowing for the incorporation of imprecise information. It also recognizes the fundamental property of coherence when updating posterior distributions, since they *can be used as the prior probability in taking into account of a further set of data*. Despite a style that is often difficult to penetrate, this is thus a major chapter of *Theory of Probability*. It will also become clearer at a later stage that the principles exposed in this chapter correspond to the (modern) notion of objective Bayes inference: despite mentions of prior probabilities as reflections of prior belief or existing pieces of information, *Theory of Probability* remains strictly "objective" in that prior distributions are always derived analytically from sampling distributions and that all examples are treated in a noninformative manner. One may find it surprising that a physicist like Jeffreys does not emphasise the appeal of subjective Bayes, that is, the ability to take into account genuine prior information in a principled way. But this is in line with both his predecessors, including Laplace and Bayes, and their use of uniform priors and his main field of study that he perceived as objective (Lindley, 2008, private communication), while one of the main appeals of *Theory of Probability* is to provide a general and coherent framework to derive objective priors.

[2]In order to keep readability as high as possible, we shall use modern notation whenever the original notation is either unclear or inconsistent, for example, Greek letters for parameters and roman letters for observations.

### 2.1 A Philosophical Exercise

The chapter starts in Section 1.0 with an epistemological discussion of the nature of (statistical) inference. Some sections are quite puzzling. For instance, the example that the kinematic equation for an object in free-fall,

$$s = a + ut + \frac{1}{2}gt^2,$$

cannot be *deduced* from observations is used as an argument against deduction under the reasoning that an infinite number of functions,

$$s = a + ut + \frac{1}{2}gt^2 + f(t)(t - t_1)\cdots(t - t_n),$$

also apply to describe a free fall observed at times $t_1, \ldots, t_n$. The limits of the epistemological discussion in those early pages are illustrated by the introduction of Ockham's razor (*the choice of the simplest law that fits the fact*), as the meaning of what a *simplest law* can be remains unclear, and the section lacks a clear (objective) argument in motivating this choice, besides common sense, while the discussion ends up with a somehow paradoxical statement that, since *deductive logic provides no explanation of the choice of the simplest law*, this is *proof that deductive logic is grossly inadequate to cover scientific and practical requirements*. On the other hand, and from a statistician's narrower perspective, one can re-interpret this gravity example as possibly the earliest discussion of the conceptual difficulties associated with model choice, which are still not entirely resolved today. In that respect, it is quite fascinating to see this discussion appear so early in the book (third page), as if Jeffreys had perceived how important this debate would become later.

Note that, maybe due to this very call to Ockham, the later Bayesian literature abounds in references to *Ockham's razor* with little formalization of this principle, even though Berger and Jefferys (1992), Balasubramanian (1997) and MacKay (2002) develop elaborate approaches. In particular, the definition of the Bayes factor in Section 1.6 can be seen as a partial implementation of Ockham's razor when setting the probabilities of both models equal to 1/2. In the beginning of his Chapter 28, entitled *Model Choice and Occam's Razor*, MacKay (2002) argues that Bayesian inference embodies Ockham's razor because "simple" models tend to produce more precise predictions and, thus, when the data is equally compatible with several models, the simplest one will end up as the most probable. This is generally true, even though there are some counterexamples in Bayesian nonparametrics.

Overall, we nonetheless feel that this part of *Theory of Probability* could be skipped at first reading as less relevant for Bayesian studies. In particular, the opposition between mathematical deduction and statistical induction does not appear to carry a strong argument, even though the distinction needs (needed?) to be made for mathematically oriented readers unfamiliar with statistics. However, from a historical point of view, this opposition must be considered against the then-ongoing debate about the nature of induction, as illustrated, for instance, by Karl Popper's articles of this period about the logical impossibility of induction (Popper, 1934).

### 2.2 Foundational Principles

The text becomes more focused when dealing with the construction of a *theory of inference*: while some notions are yet to be defined, including the pervasive *evidence*, sentences like *inference involves in its very nature the possibility that the alternative chosen as the most likely may in fact be wrong* are in line with our current interpretation of modeling and obviously with the Bayesian paradigm. In Section 1.1 Jeffreys sets up a collection of *postulates* or *rules* that act like axioms for his theory of inference, some of which require later explanations to be fully understood:

1. *All hypotheses must be explicitly stated and the conclusions must follow from the hypotheses*: what may first sound like an obvious scientific principle is in fact a leading characteristic of Bayesian statistics. While it seems to open a whole range of new

questions—"To what extent must we define our belief in the statistical models used to build our inference? How can a unique conclusion stem from a given model and a given set of observations?"—and while it may sound far too generic to be useful, we may interpret this statement as setting the working principle of Bayesian decision theory: given a prior, a sampling distribution, an observation and a loss function, there exists a single decision procedure. In contrast, the frequentist theories of Neyman or of Fisher require the choice of *ad hoc* procedures, whose (good or bad) properties they later analyze. But this may be a far-fetched interpretation of this rule at this stage even though the comment will appear more clearly later.

2. *The theory must be self-consistent.* The statement is somehow a repetition of the previous rule and it is only later (in Section 3.10) that its meaning becomes clearer, in connection with the introduction of Jeffreys's noninformative priors as a self-contained principle. Consistency is nonetheless a dominant feature of the book, as illustrated in Section 3.1 with the rejection of Haldane's prior.[3]

3. *Any rule must be applicable in practice.* This "rule" does not seem to carry any weight in practice. In addition, the explicit prohibition of estimates based on *impossible experiments* sounds implementable only through deductive arguments. But this leads to the exclusion of rules based on frequency arguments and, as such, is fundamental in setting a Bayesian framework. Alternatively (and this is another interpretation), this constraint should be worded in more formal terms of the measurability of procedures.

4. *The theory must provide explicitly for the possibility that inferences made by it may turn out to be wrong.* This is both a fundamental aspect of statistical inference and an indication of a surprising view of inference. Indeed, even when conditioning on the model, inference is never *right* in the sense that a point estimate rarely gives the *true* answer. It may be that Jeffreys is solely thinking of statistical testing, in which case the *rightfulness* of a decision is necessarily conditional on the *truthfulness* of the corresponding model and thus dubious. A more relative (or more precise) statement would have been more adequate. But, from reading further (as in Section 1.2), it appears that this *rule* is to be understood as the foundational principle (*the chief constructive rule*) for defining prior distributions. While this is certainly not clear at this stage, Bayesian inference does indeed provide for the possibility that the model under study is not correct and for the unreliability of the resulting inference via a posterior probability.

5. *The theory must not deny any empirical proposition a priori.* This principle remains unclear when put into practice. If it is to be understood in the sense of a physical theory, there is no reason why some empirical proposition could not be excluded from the start. If it is the sense of an inferential theory, then the statement would require a better definition of *empirical proposition*. But Jeffreys using the epithet a priori seems to imply that the prior distribution corresponding to the *theory* must be as inclusive as possible. This certainly makes sense as long as prior information does not exclude parts of the parameter space as, for instance, in Physics.

6. *The number of postulates should be reduced to a minimum.* This rule sounds like an embedded Ockham's razor, but, more positively, it can also be interpreted as a call for *noninformative* priors. Once again, the vagueness of the wording opens a wide range of interpretations.

7. *The theory need not represent thought-processes in details, but should agree with them in outline.* This vague principle could be an attempt at reconciliating statistical theories, but it does not give clear directions on how to proceed. In the light of Jeffreys's arguments, it could rather signify that the construction of prior distributions cannot exactly reflect an actual construction in real life. Since a noninformative (or "objective") perspective is adopted for most of the book, this is more likely to be a preliminary argument in favor of this line of thought. In Section 1.2 this rule is invoked to derive the (prior) ordering of events.

8. *An objection carries no weight if [it] would invalidate part of pure mathematics.* This rule grounds *Theory of Probability* within mathematics, which may be a necessary reminder in the spirit of the time (where some were attempting to dissociate statistics from mathematics).

The next paragraph discusses the notion of *probability*. Its interest is mostly historical: in the early 1930's, the axiomatic definition of probability based on Kolmogorov's axioms was not yet universally accepted, and there were still attempts to base this definition on limiting properties. In particular, Lebesgue integration was not part of the undergraduate curriculum till the late 1950's at either Cambridge or Oxford (Lindley, 2008, private communication). This debate is no longer relevant, and the current theory of probability, as derived from measure theory, does not bear further discussion. This also removes the ambiguity of constructing *objective probabilities* as derived from *actual or possible observations*. A probability model is to be understood as a mathematical (and thus unobjectionable) construct, in agreement with Rule 8 above.

[3]Consistency is then to be understood in the weak sense of invariant under reparameterization, which is a usual argument for Jeffreys's principle, not in terms of asymptotic convergence properties.

Then follows (still in Section 1.1) a rather long debate on *causality* versus *determinism*. While the principles stated in those pages are quite acceptable, the discussion only uses the most basic concept of *determinism*, namely, that identical causes give identical effects, in the sense of Laplace. We thus agree with Jeffreys that, at this level, *the principle is useless*, but the same paragraph actually leaves us quite confused as to its real purpose. A likely explanation (Lindley, 2008, personal communication) is that Jeffreys stresses the inevitability of probability statements in Science: (measurement) errors are not mistakes but part of the picture.

## 2.3 Prior Distributions

In Section 1.2 Jeffreys introduces the notion of prior in an indirect way, by considering that the probability of a proposition is always conditional on some data and that the occurrence of new items of information (*new evidence*) on this proposition simply updates the available data. This is slightly contrary to our current way of defining a prior distribution $\pi$ on a parameter $\theta$ as the information available on $\theta$ *prior* to the observation of the data, but it simply conveys the fact that the prior distribution must be derived from some prior items of information about $\theta$. As pointed out by Jeffreys, this also allows for the coexistence of prior distributions for different experts within the same probabilistic framework.[4] In the sequel all statements will, however, condition on the *same* data.

The following paragraphs derive standard mathematical logic axioms that directly follow from a formal (modern) definition of a probability distribution, with the provision that this probability is always conditional on the *same* data. This is also reminiscent of the derivation of the existence of a prior distribution from an ordering of prior probabilities in DeGroot (1970), but the discussion about the arbitrary ranking of probabilities between 0 and 1 may sound anecdotal today. Note also that, from a mathematical point of view, defining only conditional probabilities like $P(p|q)$ is somehow superfluous in that, if the conditioning $q$ is to remain fixed, $P(\cdot|q)$ is a regular probability distribution, while, if $q$ is to be updated into $qr$, $P(\cdot|qr)$ can be derived from $P(\cdot|q)$ by Bayes' theorem (which is to be introduced later). Therefore, in all cases, $P(\cdot|q)$ appears like the reference probability. At some stage, while stating that the probability of the sure event is equal to one is merely a convention, Jeffreys indicates that, when expressing *ignorance over an infinite range of values of a quantity, it may be convenient* to use $\infty$ instead. Clearly, this paves the way for the introduction of improper priors.[5] Unfortunately, the convention and the motivation (*to keep ratios for finite ranges determinate*) do not seem correct, if in tune with the perspective of the time (see, e.g., Lhoste, 1923; Broemeling and Broemeling, 2003). Notably, setting all events involving an infinite range with a probability equal to $\infty$ seems to restrict the abilities of the theory to a far extent.[6] Similar to Laplace, Jeffreys is more used to handling equal probability finite sets than continuous sets and the extension to continuous settings is unorthodox, using, for instance, Dedekind's sections and putting several meanings under the notation $dx$. Given the convoluted derivation of conditional probabilities in this context, the book states the product rule $P(qr|p) = P(q|p)P(r|qp)$ as an axiom, rather than as a consequence of the basic probability axioms. It leads (in Section 1.22) to Bayes' theorem,

[4] Jeffreys seems to further note that the same conditioning applies for the model of reference.

[5] Jeffreys's *Theory of Probability* strongly differs from the earlier *Scientific Inference* (1931) in this respect, the latter being rather dismissive of the mathematical difficulty: *To make this integral equal to* 1 *we should therefore have to include a zero factor unless very small and very large values are excluded. This does appear to be the case* (Section 5.43, page 67).

[6] This difficulty with handling $\sigma$-finite measures and continuous variables will be recurrent throughout the book: Jeffreys does not seem to be adverse to normalizing an improper distribution by $\infty$, even though the corresponding derivations are not meaningful.

namely, that, for all events $q_r$,

$$P(q_r|pH) \propto P(q_r|H)P(p|q_rH),$$

where $H$ denotes the *information available* and $p$ a *set of observations*. In this (modern) format $P(p|q_rH)$ is identified as Fisher likelihood and $P(q_r|H)$ as the prior probability. Bayes' theorem is defined as the *principle of inverse probability* and only for finite sets, rather than for measures.[7] Obviously, the general version of Bayes' theorem is used in the sequel for continuous parameter spaces.

Section 1.3 represents one of the few forays of the book into the realm of decision theory,[8] in connection with Laplace's notions of mathematical and moral expectations, and with Bernoulli's Saint Petersburg paradox, but there is no recognition of the central role of the loss function in defining an optimal Bayes rule as formalized later by Wald (1950) and Raiffa and Schlaifer (1961). The attribution of a decision-theoretic background to T. Bayes himself is surprising, since there is not anything close to the notion of loss or of benefit in Bayes' (1763) original paper. We nonetheless find there the seed of an idea later developed in Rubin (1987), among others, that prior and loss function are indistinguishable. [Section 1.8 briefly re-enters this perspective to point out that (posterior) expectations are often *nowhere near* the actual value of the random quantity.] The next section (Section 1.4) is important in that it tackles for the first time the issue of noninformative priors. When the number of alternatives is finite, Jeffreys picks the uniform prior as his noninformative prior, following Laplace's *Principle of Insufficient Reason*. The difficulties associated with this choice in continuous settings are not mentioned at this stage.

## 2.4 More Axiomatics and Some Asymptotics

Section 1.5 attempts an axiomatic derivation that the Bayesian principles just stated follow the rules imposed earlier. This part does not bring much novelty, once the fundamental properties of a probability distribution are stated. This is basically the purpose of this section, where earlier "Axioms" are checked in terms of the posterior probability $P(\cdot|pH)$. A reassuring consequence of this derivation is that the use of a posterior probability as the basis for inference cannot lead to inconsistency. The use of the posterior as a new prior for future observations and the corresponding learning principle are developed at this stage. The debate about the choice of the prior distribution is postponed till later, while the issue of the influence of this prior distribution is dismissed as having *very little difference [on] the results*, which needs to be quantified, as in the quote below at the beginning of Section 5.

Given the informal approach to (or rather without) measure theory adopted in *Theory of Probability*, the study of the limiting behavior of posterior distributions in Section 1.6 does not provide much insight. For instance, the fact that

$$P(q|p_1\cdots p_nH) = \frac{P(q|H)}{P(p_1|H)P(p_2|p_1H)\cdots P(p_n|p_1\cdots p_{n-1}H)}$$

is shown to induce that $P(p_n|p_1\cdots p_{n-1}H)$ converges to 1 is not particularly surprising, although it relates to Laplace's principle that *repeated verifications of consequences of a hypothesis will make it practically certain that the next consequence will be verified*. It would have been equally interesting to focus on cases in which $P(q|p_1\cdots p_nH)$ goes to 1.

The end of Section 1.62 introduces some quantities of interest, such as the distinction between estimation problems and significance tests, but with no clear guideline: when comparing models of complexity $m$ (this quantity being only defined for differential equations), Jeffreys suggests using prior probabilities that are penalized by $m$, such as $2^{-m}$ or $6/\pi^2m^2$, the motivation for those specific values being that the corresponding series converge. Penalization by the model complexity is quite an interesting idea, to be formalized later by, for example, Rissanen (1983, 1990), but Jeffreys somehow kills this idea before it is hatched by pointing out the difficulties with the definition of $m$.

Instead, Jeffreys switches to a completely different (if paramount) topic by defining in a few lines the Bayes factor for testing a point null hypothesis,

$$K = \frac{P(q|\theta H)}{P(q'|\theta H)} \bigg/ \frac{P(q|H)}{P(q'|H)},$$

[7] As noted by Fienberg (2006), the adjective term "Bayesian" had not yet appeared in the statistical literature by the time *Theory of Probability* was published, and Jeffreys sticks to the 19th century denomination of "inverse probability." The adjective can be traced back to either Ronald Fisher, who used it in a rather derogatory meaning, or to Abraham Wald, who gave it a more complimentary meaning in Wald (1950).

[8] The reference point estimator advocated by Jeffreys (if any) seems to be the maximum a posteriori (MAP) estimator, even though he stated in his discussion of Lindley (1953) that he *deprecated the whole idea of picking out a unique estimate*.

where $\theta$ denotes the data. He suggests using $P(q|H) = 1/2$ as a default value, except for sequences of embedded hypotheses for which he suggests

$$\frac{P(q|H)}{P(q'|H)} = 2,$$

presumably because the series with leading term $2^{-n}$ is converging.

Once again, the rather quick coverage of this material is somehow frustrating, as further justifications would have been necessary for the choice of the constant and so on.[9] Instead, the chapter concludes with a discussion of the distinction between "idealism" and "realism" that can be skipped for most purposes.

## 3. CHAPTER II: DIRECT PROBABILITIES

*The whole of the information contained in the observations that is relevant to the posterior probabilities of different hypotheses is summed up in the values that they give to the likelihood.*
H. JEFFREYS, *Theory of Probability*, Section 2.0.

This chapter is certainly the least "Bayesian" chapter of the book, since it covers both the standard *sampling* distributions and some equally standard probability results. It starts with a reminder that the *principle of inverse probability can be stated in the form*

$$\textit{Posterior Probability} \propto \textit{Prior Probability} \cdot \textit{Likelihood},$$

thus rephrasing Bayes' theorem in terms of the likelihood and with the proper indication that the *relevant information contained in the observations* is summarized by the likelihood (*sufficiency* will be mentioned later in Section 3.7). Then follows (still in Section 2.0) a long paragraph about the tentative nature of models, concluding that a statistical model must be made part of the prior information $H$ before it can be tested against the observations, which (presumably) relates to the fact that Bayesian model assessment must involve a description of the alternative(s) to be validated.

The main bulk of the chapter is about sampling distributions. Section 2.1 introduces binomial and hypergeometric distributions at length, including the interesting problem of deciding between binomial versus negative binomial experiments when faced with the outcome of a survey, used later in the defence of the Likelihood Principle (Berger and Wolpert, 1988). The description of the binomial contains the equally interesting remark that a given coin repeatedly thrown will show a bias toward head or tail due to the wear, a remark later exploited in Diaconis and Ylvisaker (1985) to justify the use of mixtures of conjugate priors. Bernoulli's version of the Central Limit theorem is also recalled in this section, with no particular appeal if one considers that a modern Statistics course (see, e.g., Casella and Berger, 2001) would first start with the probabilistic background.[10]

The Poisson distribution is first introduced as a limiting distribution for the binomial distribution $\mathcal{B}(n, p)$ when $n$ is large and $np$ is bounded. (Connections with radioactive disintegration are mentioned afterward.) The normal distribution is proposed as a large sample approximation to a sum of Bernoulli random variables. As for the other distributions, there is some attempt at justifying the use of the normal distribution, as well as [what we find to be] a confusing paragraph about the "true" and "actual observed" values of the parameters. A long section (Section 2.3) expands about the properties of Pearson's distributions, then allowing Jeffreys to introduce the negative binomial as a mixture of Poisson distributions. The introduction of the bivariate normal distribution is similarly convoluted, using first binomial variates and second a limiting argument, and without resorting to matrix formalism.

Section 2.6 attempts to introduce cumulative distribution functions in a more formal manner, using the current three-step definition, but again dealing with limits in an informal way. Rather coherently from a geophysicist's point of view, characteristic functions are also covered in great detail, including connections with moments and the Cauchy distribution, as well as Lévy's inversion theorem. The main goal of using characteristic functions seems nonetheless to be able to establish the Central Limit theorem in its full generality (Section 2.664).

[9]Similarly, the argument against *philosophers that maintain that no method based on the theory of probability can give a (…) non-zero probability to a precise value against a continuous background* is not convincing as stated. The distinction between zero measure events and mixture priors including a Dirac mass should have been better explained, since this is the basis for Bayesian point-null testing.

[10]In fact, some of the statements in *Theory of Probability* that surround the statement of the Central Limit theorem are not in agreement with measure theory, as, for instance, the confusion between pointwise and uniform convergence, and convergence in probability and convergence in distribution.

Rather surprisingly for a Bayesian reference book and mostly in complete disconnection with the testing chapters, the $\chi^2$ test of goodness of fit is given a large and uncritical place within this book, including an adjustment for the degrees of freedom.[11] Examples include the obvious independence of a rectangular contingency table. The only criticism (Section 2.76) is fairly obscure in that it blames poor performances of the $\chi^2$ test on the fact that all divergences in the $\chi^2$ sum are equally weighted. The test is nonetheless implemented in the most classical manner, namely, that the hypothesis is rejected if the $\chi^2$ statistic is outside the standard interval. It is unclear from the text in Section 2.76 that rejection would occur were the $\chi^2$ statistic too small, even though Jeffreys rightly addresses the issue at the end of Chapter 5 (Section 5.63). He also mentions the need to coalesce small groups into groups of size at least 5 with no further justification. The chapter concludes with similar uses of Student's $t$ and Fisher's $z$ tests.

## 4. CHAPTER III: ESTIMATION PROBLEMS

*If we have no information relevant to the actual value of the parameter, the probability must be chosen so as to express the fact that we have none.*
H. Jeffreys, *Theory of Probability*, Section 3.1.

This is a major chapter of *Theory of Probability* as it introduces both exponential families and the principle of Jeffreys noninformative priors. The main concepts are already present in the early sections, including some invariance principles. The purpose of the chapter is stated as a point estimation problem, where obtaining the *probability distribution of [the] parameters, given the observations* is the goal. Note that *estimation* is not to be understood in the (modern?) sense of point estimation, that is, as a way to produce numerical substitutes for the true parameters that are based on the data, since the decision-theoretic perspective for building (point) estimators is mostly missing from the book (see Section 1.8 for a very brief remark on expectations). Both Good (1980) and Lindley (1980) stress this absence.

### 4.1 Noninformative Priors of Former Days

Section 3.1 sets the principles for selecting noninformative priors. Jeffreys recalls Laplace's rule that, if a parameter is real-valued, its *prior probability should be taken as uniformly distributed*, while, if this parameter is positive, the *prior probability of its logarithm should be taken as uniformly distributed*. The motivation advanced for using both priors is the *invariance principle*, namely, the invariance of the prior selection under *several different sets of parameters*. At this stage, there is no recognition of a potential problem with using a $\sigma$-finite measure and, in particular, with the fact that these priors are not probability distributions, but rather a simple warning that these are *formal rules expressing ignorance*. We face the difficulty mentioned earlier when considering $\sigma$-finite measures since they are not properly handled at this stage: when stating that one starts *with any distribution of prior probability*, it is not possible to include $\sigma$-finite measures this way, except via the [incorrect] argument that *a probability is merely a number* and, thus, that the total weight can be $\infty$ as well as 1: *use* $\infty$ *instead of* 1 *to indicate certainty on data* $H$. The wrong interpretation of a $\sigma$-finite measure as a probability distribution (and of $\infty$ as a "number") then leads to immediate paradoxes, such as the prior probability of any finite range being null, which sounds *inconsistent with the statement that we know nothing about* the parameter, but this results from an over-interpretation of the measure as a probability distribution already pointed out by Lindley (1971, 1980) and Kass and Wasserman (1996).

The argument for using a flat (Lebesgue) prior is based (a) on its use by both Bayes and Laplace in finite or compact settings, and (b) on the argument that it correctly reflects the absence of prior knowledge *about the value of the parameter*. At this stage, no point is made against it for reasons related with the *invariance principle*—there is only one parameterization that coincides with a uniform prior—but Jeffreys already argues that flat priors cannot be used for significance tests, because they would always reject the point null hypothesis. Even though Bayesian significance tests, including Bayes factors, have not yet been properly introduced, the notion of an infinite mass canceling a point null hypothesis is sufficiently intuitive to be used at this point.

[11] Interestingly enough, the parameters are estimated by minimum $\chi^2$ rather than either maximum likelihood or Bayesian point estimates. This is, again, a reflection of the practice of the time, coupled with the fact that most approaches are asymptotically indistinguishable. Posterior expectations are not at all advocated as Bayes (point) estimators in *Theory of Probability*.

While, indeed, using an improper prior is a major difficulty when testing point null hypotheses because it gives an infinite mass to the alternative (DeGroot, 1970), Jeffreys fails to identify the problem as such but rather blames the flat prior *applied to a parameter with a semi-infinite range of possible values*. He then goes on justifying the use of $\pi(\sigma) = 1/\sigma$ for positive parameters (replicating the argument of Lhoste, 1923) on the basis that it is invariant for the change of parameters $\varrho = 1/\sigma$, as well as any other power, failing to recognize that other transforms that preserve positivity do not exhibit such an invariance. One has to admit, however, that, from a physicist's perspective, power transforms are more important than other mathematical transforms, such as arctan, because they can be assigned meaningful units of measurement, while other functions cannot. At least this seems to be the spirit of the examples considered in *Theory of Probability*: *Some methods of measuring the charge of an electron give e, others $e^2$*.

There is a vague indication that Jeffreys may also recognize $\pi(\sigma) = 1/\sigma$ as the scale group invariant measure, but this is unclear. An indefensible argument follows, namely, that

$$\int_0^a v^n \, dv \Big/ \int_a^\infty v^n \, dv$$

is only indeterminate when $n = -1$, which allows us to avoid contradictions about the lack of prior information. Jeffreys acknowledges that this does not solve the problem since this choice implies that the prior "probability" of a finite interval $(a, b)$ is then always null, but he avoids the difficulty by admitting that the probability that $\sigma$ falls in a particular range is zero, because *zero probability does not imply impossibility*. He also acknowledges that the *invariance principle* cannot encompass the whole range of transforms without being *inconsistent*, but he nonetheless sticks to the $\pi(\sigma) = 1/\sigma$ prior as it is *better than the Bayes–Laplace rule*.[12] Once again, the argument sustaining the whole of Section 3.1 is incomplete since missing the fundamental issue of distinguishing proper from improper priors.

While Haldane's (1932) prior on probabilities (or rather on *chances* as defined in Section 1.7),

$$\pi(p) \propto \frac{1}{p(1-p)},$$

is dismissed as too extreme (and *inconsistent*), there is no discussion of the main difficulty with this prior (or with any other improper prior associated with a finite-support sampling distribution), which is that the corresponding posterior distribution is not defined when $x \sim \mathcal{B}(n, p)$ is either equal to 0 or to $n$ (although Jeffreys concludes that $x = 0$ leads to a point mass at $p = 0$, due to the infinite mass normalization).[13] Instead, the corresponding Jeffreys's prior

$$\pi(p) \propto \frac{1}{\sqrt{p(1-p)}}$$

is suggested with little justification against the (truly) uniform prior: *we may as well use the uniform distribution*.

## 4.2 Laplace's Succession Rule

Section 3.2 contains a Bayesian processing of Laplace's succession rule, which is an easy introduction given that the parameter of the sampling distribution, a hypergeometric $\mathcal{H}(N, r)$, is an integer. The choice of a uniform prior on $r$, $\pi(r) = 1/(N+1)$, does not require much of a discussion and the posterior distribution

$$\pi(r|l, m, N, H) = \binom{r}{l}\binom{N-r}{m} \Big/ \binom{N+1}{l+m+1}$$

is available in closed form, including the normalizing constant. The posterior predictive probability that *the next specimen will be of the same type* is then $(l+1)/(l+m+1)$ and more complex predictive probabilities can be computed as well. As in earlier books involving Laplace's succession rule, the section argues about its truthfulness from a metaphysical point of view (using classical arguments about the probabilities that *the sun rising tomorrow* and that *all swans are white* that always seem to be associates themselves with this topic) but, more interestingly, it then moves to introducing a point mass on specific values of the parameter in preparation for hypothesis testing. Namely, following a renewed criticism of the *uniform assessment* via the fact that

$$\frac{P(r=N|l,m=0,N,H)}{P(r\neq N|l=n,N,H)}=\frac{l+1}{N+1}$$

is too small, Jeffreys suggests setting aside a portion $2k$ of the prior mass for both extreme values $r=0$ and $r=N$. This is indeed equivalent to using a point mass on the null hypothesis of homogeneity of the population. While mixed samples are independent of the choice of $k$ (since they exclude those extreme values), a sample *of the first type* with $l=n$ leads to a posterior probability ratio of

$$\frac{P(r=N|l=n,N,H)}{P(r\neq N|l=n,N,H)}=\frac{n+1}{N-n}\frac{k}{1-2k}\frac{N-1}{1},$$

which leads to the crucial question of the choice[14] of $k$. The ensuing discussion is not entirely convincing: $\frac{1}{2}$ *is too large*, $\frac{1}{4}$ *is not unreasonable [but] too low in this case*. The alternative

$$k=\frac{1}{4}+\frac{1}{N+1}$$

argues that the *classification of possibilities [is] as follows*: (1) *Population homogeneous on account of some general rule.* (2) *No general rule but extreme values to be treated on a level with others.* This proposal is mostly interesting for its bearing on the continuous case, for, in the finite case, it does not sound logical to put weight on the null hypothesis ($r=0$ and $r=N$) within the alternative, since this confuses the issue. (See Berger, Bernardo and Sun, 2009, for a recent reappraisal of this approach from the point of view of reference priors.)

[12] In both the 19th and early 20th centuries, there is a tradition within the not-yet-Bayesian literature to go to extreme lengths in the justification of a particular prior distribution, as if there existed one golden prior. See, for example, Broemeling and Broemeling (2003) in this respect.

[13] Jeffreys (1931, 1937) does address the problem in a clearer manner, stating that *this is not serious, for so long as the sample is homogeneous (meaning $x = 0, n$) the extreme values (meaning $p = 0, 1$) are still admissible, and we do attach a high probability to the proposition is of one type; while as soon as any exceptions are known the extreme values are completely excluded and no infinity arises* (Section 10.1, page 195).

Section 3.3 seems to extend Laplace's succession rule to the case in which *the class sampled consists of several types*, but it actually deals with the (much more interesting) case of Bayesian inference for the multinomial $\mathcal{M}(n;p_1,\ldots,p_r)$ distribution, when using the Dirichlet $\mathcal{D}(1,\ldots,1)$ distribution as a prior. Jeffreys recovers the Dirichlet $\mathcal{D}(x_1+1,\ldots,x_r+1)$ distribution as the posterior distribution and he derives the predictive probability that *the next member will be of the first type* as

$$(x_1+1)/\sum_i x_i+r.$$

There could be some connections there with the irrelevance of alternative hypotheses later (in time) discussed in polytomous regression models (Gouriéroux and Monfort, 1996), but they are well hidden. In any case, the Dirichlet distribution is not invariant to the introduction of new types.

### 4.3 Poisson Distribution

The processing of the *estimation of the parameter* $\alpha$ of the Poisson distribution $\mathcal{P}(\alpha)$ is based on the [improper] prior $\pi(\alpha)\propto 1/\alpha$, deemed to be *the correct prior probability distribution* for scale invariance reasons. Given $n$ observations from $\mathcal{P}(\alpha)$ with sum $S_n$, Jeffreys reproduces Haldane's (1932) derivation of the Gamma posterior $\mathcal{G}a(S_n,n)$ and he notes that $S_n$ is a sufficient statistic, but does not make a general property of it at this stage. (This is done in Section 3.7.)

The alternative choice $\pi(\alpha)\propto 1/\sqrt{\alpha}$ will be later justified in Section 3.10 not as Jeffreys's (invariant) prior but as leading to a posterior defined for all observations, which is not the case of $\pi(\alpha)\propto 1/\alpha$ when $x=0$, a fact overlooked by Jeffreys. Note that $\pi(\alpha)\propto 1/\alpha$ can nonetheless be advocated by Jeffreys on the ground that the Poisson process derives from the exponential distribution, for which $\alpha$ is a scale parameter: $e^{-\alpha t}$ *represents the fraction of the atoms originally present that survive after time* $t$.

### 4.4 Normal Distribution

When the sampling variance $\sigma^2$ of a normal model $\mathcal{N}(\mu,\sigma^2)$ is known, the posterior distribution associated with a flat prior is correctly derived as $\mu|x_1,\ldots,x_n\sim\mathcal{N}(\bar{x},\sigma^2/n)$ (with the repeated difficulty about the use of a $\sigma$-finite measure as a probability). Under the joint improper prior

$$\pi(\mu,\sigma)\propto 1/\sigma,$$

the (marginal) posterior on $\mu$ is obtained as a Student's $t$

$$\mathcal{T}(n-1,\bar{x},s^2/n(n-1))$$

distribution, while the marginal posterior on $\sigma^2$ is an inverse gamma $\mathcal{IG}((n-1)/2,s^2/2)$.[15]

[14] A prior weight of $2k=1/2$ is reasonable since it gives equal probability to both hypotheses.

[15] Section 3.41 also contains the interesting remark that, conditional on two observations, $x_1$ and $x_2$, the posterior

Jeffreys notices that, when $n=1$, the above prior does not lead to a proper posterior since $\pi(\mu|x_1) \propto 1/|\mu - x_1|$ is not integrable, but he concludes that *the solution degenerates in the right way,* which, we suppose, is meant to say that there is not enough information in the data. But, without further formalization, it is a delicate conclusion to make.

Under the same noninformative prior, the predictive density of a second sample with sufficient statistic $(\bar{x}_2, s_2)$ is found[16] to be proportional to

$$\left\{n_1 s_1^2 + n_2 s_2^2 + \frac{n_1 n_2}{n_1+n_2}(\bar{x}_2 - \bar{x}_1)^2\right\}^{-(n_1+n_2-1)/2}.$$

A direct conclusion is that this implies that $\bar{x}_2$ and $s_2$ are dependent for the predictive, if *independent given* $\mu$ *and* $\sigma$, while the marginal predictives on $\bar{x}_2$ and $s_2^2$ are Student's $t$ and Fisher's $z$, respectively. Extensions to the prediction of multiple future samples with the same (Section 3.43) or with different (Section 3.44) means follow without surprise. In the latter case, given $m$ samples of $n_r$ $(1 \le r \le m)$ normal $\mathcal{N}(\mu_i, \sigma^2)$ measurements, the posterior on $\sigma^2$ under the noninformative prior

$$\pi(\mu_1, \ldots, \mu_r, \sigma) \propto 1/\sigma$$

is again an inverse gamma $\mathcal{IG}(\nu/2, s^2/2)$ distribution,[17] with $s^2 = \sum_r \sum_i (x_{ri} - \bar{x}_r)^2$ and $\nu = \sum_r n_r$, while the posterior on $t = \sqrt{n_i}(\mu_i - \bar{x}_i)/s$ is a Student's $t$ with $\nu$ degrees of freedom for all $i$'s (no matter what the number of observations within this group is). Figure 1 represents the posteriors on the means $\mu_i$ for the data set analyzed in this section on seven sets of measurements of the gravity. A paragraph in Section 3.44 contains hints about hierarchical Bayes modeling as a way of strengthening estimation, which is a perspective later advanced in favor of this approach (Lindley and Smith, 1972; Berger and Robert, 1990).

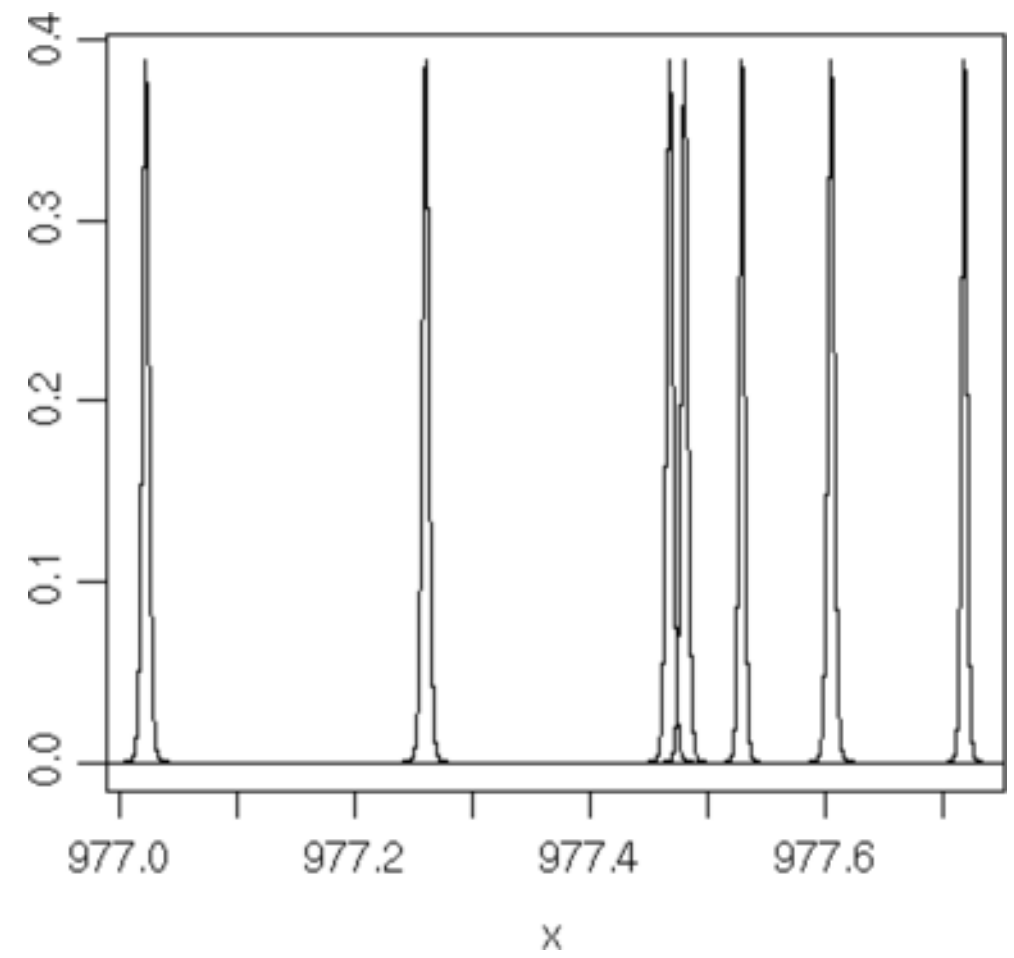

FIG. 1. *Seven posterior distributions on the values of acceleration due to gravity (in cm/sec$^2$) at locations in East Africa when using a noninformative prior.*

The extension in Section 3.5 to the setting of the normal linear regression model should be simple (see, e.g., Marin and Robert, 2007, Chapter 3), except that the use of tensorial conventions—like *when a suffix i is repeated it is to be given all values from* 1 *to m*—and the absence of matrix notation makes the reading quite arduous for today's readers.[18] Because of this lack of matrix tools, Jeffreys uses an implicit diagonalization of the regressor matrix $X^{\mathrm{T}}X$ (with modern notation) and thus expresses the posterior in terms of the transforms $\xi_i$ of the regression coefficients $\beta_i$. This section is worth reading if only to realize the immense advantage of using matrix notation. The case of regression equations

$$y_i = X_i\beta + \varepsilon_i, \quad \varepsilon_i \sim \mathcal{N}(0, \sigma_i^2),$$

with different unknown variances leads to a poly-$t$ output (Bauwens, 1984) under a noninformative prior, which is deemed to be a *complication*, and Jeffreys prefers to revert to the case when $\sigma_i^2 = \omega_i \sigma^2$ with known $\omega_i$'s.[19] The final part of this section probability that $\mu$ is between both observations is exactly 1/2. Jeffreys attributes this property to the fact that the scale $\sigma$ is directly estimated from those two observations under a noninformative prior. Section 3.8 generalizes the observation to all location-scale families with median equal to the location. Otherwise, the posterior probability is less than 1/2. Similarly, the probability that a third observation $x_3$ will be between $x_1$ and $x_2$ is equal to 1/3 under the predictive. While Jeffreys gives a proof by complete integration, this is a direct consequence of the exchangeability of $x_1$, $x_2$ and $x_3$. Note also that this is one of the rare occurrences of a credible interval in the book.

[16]In the current 1961 edition, $n_2 s_2^2$ is mistakenly typed as $n_2^2 s_2^2$ in equation (6) of Section 3.42.

[17]Jeffreys does not use the term "inverse gamma distribution" but simply notes that this is a distribution with a scale parameter that *is given by a single set of tables* (for a given $\nu$). He also notices that the distribution of the transform $\log(\sigma/s)$ is closer to a normal distribution than the original.

[18]Using the notation $c_i$ for $y_i$, $x_i$ for $\beta_i$, $y_i$ for $\hat{\beta}_i$ and $a_{ir}$ for $x_{ir}$ certainly makes reading this part more arduous.

[19]Sections 3.53 and 3.54 detail the numerical resolution of the normal equations by iterative methods and have no real bearing on modern Bayesian analysis.

mentions the interesting subcase of estimating a normal mean $\alpha$ when *truncated at* $\alpha = 0$: negative observations do not need to be rejected since only the posterior distribution has to be truncated in 0. [In a similar spirit, Section 3.6 shows how to process a uniform $\mathcal{U}(\alpha-\sigma, \alpha+\sigma)$ distribution under the noninformative $\pi(\alpha, \sigma) = 1/\sigma$ prior.]

Section 3.9 examines the estimation of a two-dimensional covariance matrix

$$\Theta = \begin{pmatrix} \sigma^2 & \varrho\sigma\tau \\ \varrho\sigma\tau & \tau^2 \end{pmatrix}$$

under centred normal observations. The prior advocated by Jeffreys is $\pi(\tau, \sigma, \varrho) \propto 1/\tau\sigma$, leading to the (marginal) posterior

$$\begin{aligned}
&\pi(\varrho|\hat{\varrho}, n) \\
&\quad \propto \int_0^\infty \frac{(1-\varrho^2)^{n/2}}{(\cosh\beta - \varrho\hat{\varrho})^n}\, d\beta \\
&\quad = \frac{(1-\varrho^2)^{n/2}}{(1-\varrho\hat{\varrho})^{n-1/2}} \\
&\qquad \cdot \int_0^1 \frac{(1-u)^{n-1}}{\sqrt{2u}} \{1-(1+\varrho\hat{\varrho})u/2\}^{-1/2}\, du
\end{aligned}$$

that only depends on $\hat{\varrho}$. (Jeffreys notes that, when $\sigma$ and $\tau$ are known, the posterior of $\varrho$ also depends on the empirical variances for both components. This paradoxical increase in the dimension of the sufficient statistics when the number of parameters is decreasing is another illustration of the limited meaning of marginal sufficient statistics pointed out by Basu, 1988.) While this integral can be computed via confluent hypergeometric functions (Gradshteyn and Ryzhik, 1980),

$$\begin{aligned}
&\int_0^1 \frac{(1-x)^{n-1}}{\sqrt{u(1-au)}}\, du \\
&\quad = B(1/2, n)\, {}_2F_1\{1/2, 1/2; n+1/2; (1+\varrho\hat{\varrho})/2\},
\end{aligned}$$

the corresponding posterior is certainly less manageable than the inverse Wishart that would result from a power prior $|\Theta|^\gamma$ on the matrix $\Theta$ itself. The extension to noncentred observations with flat priors on the means induces a small change in the outcome in that

$$\begin{aligned}
\pi(\varrho|\hat{\varrho}, n) &\propto \frac{(1-\varrho^2)^{(n-1)/2}}{(1-\varrho\hat{\varrho})^{n-3/2}} \\
&\quad \cdot \int_0^1 \frac{(1-u)^{n-2}}{\sqrt{2u}} \\
&\qquad \cdot \{1-(1+\varrho\hat{\varrho})u/2\}^{-1/2}\, du,
\end{aligned}$$

which is also the posterior obtained directly from the distribution of $\hat{\varrho}$. Indeed, the sampling distribution is given by

$$\begin{aligned}
f(\hat{\varrho}|\varrho) &= \frac{n-2}{\sqrt{2\pi}} (1-\hat{\varrho}^2)^{(n-4)/2} \\
&\quad \cdot (1-\varrho^2)^{(n-1)/2} \frac{\Gamma(n-1)}{\Gamma(n-1/2)} \\
&\quad \cdot (1-\varrho\hat{\varrho})^{-(n-3/2)} \\
&\quad \cdot {}_2F_1\{1/2, 1/2; n-1/2; (1+\varrho\hat{\varrho})/2\}.
\end{aligned}$$

There is thus no marginalization paradox (Dawid, Stone and Zidek, 1973) for this prior selection, while one occurs for the alternative choice $\pi(\tau, \sigma, \varrho) \propto 1/\tau^2\sigma^2$.

### 4.5 Sufficiency and Exponential Families

Section 3.7 generalizes[20] observations made previously about sufficient statistics for particular distributions (Poisson, multinomial, normal, uniform). If there exists a sufficient statistic $T(x)$ when $x \sim f(x|\alpha)$, the posterior distribution on $\alpha$ only depends on $T(x)$ and on the number $n$ of observations.[21] The generic form of densities from exponential families

$$\log f(x|\alpha) = (x-\alpha)\mu'(\alpha) + \mu(\alpha) + \psi(x)$$

is obtained by a convoluted argument of imposing $\bar{x}$ as the MLE of $\alpha$, which is not equivalent to requiring $\bar{x}$ to be sufficient. The more general formula

$$\begin{aligned}
&f(x|\alpha_1, \ldots, \alpha_m) \\
&\quad = \phi(\alpha_1, \ldots, \alpha_m)\psi(x) \exp \sum_{s=1}^m u_s(\alpha) v_s(x)
\end{aligned}$$

is provided as a consequence of the (then very recent) Pitman–Koopman[–Darmois] theorem[22] on the necessary and sufficient connection between the existence of fixed dimensional sufficient statistics and exponential families. The theorem as stated does not impose a fixed support on the densities $f(x|\alpha)$ and this invalidates the necessary part, as shown in Section 3.6 with the uniform distribution. It is only later in Section 3.6 that parameter-dependent supports are mentioned, with an unclear conclusion. Surprisingly, this section does not contain any indication that the specific structure of exponential families could be used to construct conjugate[23] priors (Raiffa, 1968). This lack of connection with regular priors highlights the fully noninformative perspective advocated in *Theory of Probability*, despite comments (within the book) that priors should reflect prior beliefs and/or information.

[20] Jeffreys's derivation remains restricted to the unidimensional case.

[21] Stating that *$n$ is an ancillary statistic* is both formally correct in Fisher's sense ($n$ does not depend on $\alpha$) and ambiguous from a Bayesian perspective since the posterior on $\alpha$ depends on $n$.

[22] Darmois (1935) published a version (in French) of this theorem in 1935, about a year before both Pitman (1936) and Koopman (1936).

### 4.6 Predictive Densities

Section 3.8 contains the rather amusing and not well-known result that, for any location-scale parametric family such that the location parameter is the median, the posterior probability that the third observation lies between the first two observations is 1/2. This may be the first use of Bayesian predictive distributions, that is, $p(x_3|x_1, x_2)$ in this case, where parameters are integrated out. Such predictive distributions cannot be properly defined in frequentist terms; at best, one may take $p(x_3|\theta = \hat{\theta})$ where $\hat{\theta}$ is a plug-in estimator. Building more sensible predictives seems to be one major appeal of the Bayesian approach for modern practitioners, in particular, econometricians.

### 4.7 Jeffreys's Priors

Section 3.10 introduces Fisher information as a quadratic approximation to distributional distances. Given the Hellinger distance and the Kullback–Leibler divergence,

$$d_1(P, P') = \int |(dP)^{1/2} - (dP')^{1/2}|^2$$

and

$$d_2(P, P') = \int \log \frac{dP}{dP'} d(P - P'),$$

we have the second-order approximations

$$d_1(P_\alpha, P_{\alpha'}) \approx \frac{1}{4}(\alpha - \alpha')^{\mathrm{T}} I(\alpha)(\alpha - \alpha')$$

and

$$d_2(P_\alpha, P_{\alpha'}) \approx (\alpha - \alpha')^{\mathrm{T}} I(\alpha)(\alpha - \alpha'),$$

where

$$I(\alpha) = \mathbb{E}_\alpha\left[\frac{\partial f(x|\alpha)}{\partial \alpha}\frac{\partial f(x|\alpha)^{\mathrm{T}}}{\partial \alpha}\right]$$

is Fisher information.[24] A first comment of importance is that $I(\alpha)$ is equivariant under reparameterization, because both distances are functional distances and thus *invariant for all nonsingular transformations of the parameters*. Therefore, if $\alpha'$ is a (differentiable) transform of $\alpha$,

$$I(\alpha') = \frac{d\alpha}{d\alpha'} I(\alpha) \frac{d\alpha^{\mathrm{T}}}{d\alpha'},$$

and this is the spot where Jeffreys states his general principle for deriving noninformative priors (Jeffreys's priors):[25]

$$\pi(\alpha) \propto |I(\alpha)|^{1/2}$$

is thus an ideal prior in that it is invariant under any (differentiable) transformation.

Quite curiously, there is no motivation for this choice of priors other than invariance (at least at this stage) and consistency (at the end of the chapter). Fisher information is only perceived as a second order approximation to two functional distances, with no connection with either the curvature of the likelihood or the variance of the score function, and no mention of the information content at the current value of the parameter or of the local discriminating power of the data. Finally, no connection is made at this stage with Laplace's approximation (see Section 4.0). The motivation for centering the choice of the prior at $I(\alpha)$ is thus uncertain. No mention is made either of the potential use of those functional distances as intrinsic loss functions for the [point] estimation of the parameters (Le Cam, 1986; Robert, 1996). However, the use of these intrinsic divergences (measures of discrepancy) to introduce $I(\alpha)$ as a key quantity seems to indicate that Jeffreys understood $I(\alpha)$ as a local discriminating power of the model and to some extent as the intrinsic factor used to compensate for the lack of invariance of $|\alpha - \alpha'|^2$. It corroborates the fact that Jeffreys's priors are known to behave particularly well in one-dimensional cases.

Immediately, a problem associated with this generic principle is spotted by Jeffreys for the normal distribution $\mathcal{N}(\mu, \sigma^2)$. While, when considering $\mu$ and $\sigma$

[23] As pointed to us by Dennis Lindley, Section 1.7 comes close to the concept of exchangeability when introducing *chances*.

[24] Jeffreys uses an infinitesimal approximation to derive $I(\alpha)$ in *Theory of Probability*, which is thus not defined this way, nor connected with Fisher.

[25] Obviously, those priors are not called *Jeffreys's priors* in the book but, as a counter-example to Steve Stigler's *law of eponymy* (Stigler, 1999), the name is now correctly associated with the author of this new concept.

separately, one recovers the invariance priors $\pi(\mu) \propto 1$ and $\pi(\sigma) \propto 1/\sigma$, Jeffreys's prior on the pair $(\mu, \sigma)$ is $\pi(\mu, \sigma) \propto 1/\sigma^2$. If, instead, $m$ normal observations with the same variance $\sigma^2$ were proposed, they would lead to $\pi(\mu_1, \ldots, \mu_m, \sigma) \propto 1/\sigma^{m+1}$, which *is unacceptable* (because it induces a growing departure from the true value as $m$ increases). Indeed, if one considers the likelihood

$$L(\mu_1, \ldots, \mu_m, \sigma) \propto \sigma^{-mn} \exp -\frac{n}{2\sigma^2} \sum_{i=1}^m \{(\bar{x}_i - \mu_i)^2 + s_i^2\},$$

the marginal posterior on $\sigma$ is

$$\sigma^{-mn-1} \exp -\frac{n}{2\sigma^2} \sum_{i=1}^m s_i^2,$$

that is,

$$\sigma^{-2} \sim \mathcal{G}a\left\{(mn-1)/2, n \sum_i s_i^2/2\right\}$$

and

$$\mathbb{E}[\sigma^2] = \frac{n \sum_{i=1}^m s_i^2}{mn-1}$$

whose own expectation is

$$\frac{mn-m}{mn-1}\sigma_0^2,$$

if $\sigma_0$ denotes the "true" standard deviation. If $n$ is small against $m$, the bias resulting from this choice will be important.[26] Therefore, in this special case, Jeffreys proposes a *departure from the general rule* by using $\pi(\mu, \sigma) \propto 1/\sigma$. (There is a further mention of difficulties with a large number of parameters when using one single scale parameter, with the same solution proposed. There may even be an indication about reference priors at this stage, when stating that some transforms do not need to be considered.)

The arc-sine law on probabilities,

$$\pi(p) = \frac{1}{\pi} \frac{1}{\sqrt{p(1-p)}},$$

is found to be the corresponding reference distribution, with a more severe criticism of the other distributions (see Section 4.1): *both the usual rule and Haldane's rule are rather unsatisfactory*. The corresponding Dirichlet $\mathcal{D}(1/2, \ldots, 1/2)$ prior is obtained on the probabilities of a multinomial distribution. Interestingly too, Jeffreys derives most of his priors by recomputing the $L_2$ or Kullback distance and by using a second-order approximation, rather than by following the genuine definition of the Fisher information matrix. Because Jeffreys's prior on the Poisson $\mathcal{P}(\lambda)$ parameter is $\pi(\lambda) \propto 1/\sqrt{\lambda}$, there is some attempt at justification, with the mention that *general rules for the prior probability give a starting point*, that is, act like reference priors (Berger and Bernardo, 1992).

In the case of the (normal) correlation coefficient, the posterior corresponding to Jeffreys's prior $\pi(\varrho, \tau, \sigma) \propto 1/\tau\sigma(1-\varrho^2)^{3/2}$ is not properly defined for a single observation, but Jeffreys does not expand on the generic improper nature of those prior distributions. In an attempt close to defining a reference prior, he notices that, with both $\tau$ and $\sigma$ fixed, the (conditional) prior is

$$\pi(\varrho) \propto \frac{\sqrt{1+\varrho^2}}{1-\varrho^2},$$

which, while improper, can also be compared to the arc-sine prior

$$\pi(\varrho) = \frac{1}{\pi} \frac{1}{\sqrt{1-\varrho^2}},$$

which is integrable as is. Note that Jeffreys does not conclude in favor of one of those priors: *We cannot really say that any of these rules is better than the uniform distribution*.

In the case of exponential families with natural parameter $\beta$,

$$f(x|\beta) = \psi(x)\phi(\beta) \exp \beta v(x),$$

Jeffreys does not take advantage of the fact that Fisher information is available as a transform of $\phi$, indeed,

$$I(\beta) = \partial^2 \log \phi(\beta)/\partial\beta^2,$$

but rather insists on the invariance of the distribution under location-scale transforms, $\beta = k\beta' + l$, which does not correctly account for potential boundaries on $\beta$.

Somehow, surprisingly, rather than resorting to the natural "Jeffreys's prior," $\pi(\beta) \propto |\partial^2 \log \phi(\beta)/$

[26]As pointed out to us by Lindley (2008, private communication), Jeffreys expresses more clearly the difficulty that *the corresponding* $t$ *distribution would always be [of index]* $(n+1)/2$*, no matter how many true values were estimated*, that is, that the natural reduction of the degrees of freedom with the number of nuisance parameters does not occur with this prior.

$\partial\beta^2|^{1/2}$, Jeffreys prefers to use the "standard" flat, log-flat and symmetric priors depending on the range of $\beta$. He then goes on to study the alternative of defining the noninformative prior via the mean parameterization suggested by Huzurbazar (see Huzurbazar, 1976),

$$\mu(\beta) = \int v(x) f(x|\beta)\,dx.$$

Given the overall invariance of Jeffreys's priors, this should not make any difference, but Jeffreys chooses to pick priors depending on the range of $\mu(\beta)$. For instance, this leads him once again to promote the Dirichlet $\mathcal{D}(1/2, 1/2)$ prior on the probability $p$ of a binomial model if considering that $\log p/(1-p)$ is unbounded,[27] and the uniform prior if considering that $\mu(p) = np$ varies on $(0, \infty)$. It is interesting to see that, rather than sticking to a generic principle inspired by the Fisher information that Jeffreys himself recognizes as *consistent* and that offers an almost universal range of applications, he resorts to group invariant (Haar) measures when *the rule, though consistent, leads to results that appear to differ too much from current practice.*

We conclude with a delicate example that is found within Section 3.10. Our interpretation of *a set of quantitative laws* $\phi_r$ *with chances* $\alpha_r$ *[such that] if* $\phi_r$ *is true, the chance of a variable* $x$ *being in a range* $dx$ *is* $f_r(x, \alpha_{r1}, \ldots, \alpha_{rn})\,dx$ is that of a mixture of distributions,

$$x \sim \sum_{r=1}^{m} \alpha_r f_r(x, \alpha_{r1}, \ldots, \alpha_{rn}).$$

Because of the complex shape (convex combination) of the distribution, the Fisher information is not readily available and Jeffreys suggests assigning a reference prior to the weights $(\alpha_1, \ldots, \alpha_m)$, that is, a Dirichlet $\mathcal{D}(1/2, \ldots, 1/2)$, along with separate reference priors on the $\alpha_{rs}$. Unfortunately, this leads to an improper posterior density (which integrates to infinity). In fact, mixture models do not allow for independent improper priors on their components (Marin, Mengersen and Robert, 2005).

## 5. CHAPTER IV: APPROXIMATE METHODS AND SIMPLIFICATIONS

*The difference made by any ordinary change of the prior probability is comparable with the effect of one extra observation.*

H. Jeffreys, *Theory of Probability*, Section 4.0.

As in Chapter II, many points of this chapter are outdated by modern Bayesian practice. The main bulk of the discussion is about various approximations to (then) intractable quantities or posteriors, approximations that have limited appeal nowadays when compared with state-of-the-art computational tools. For instance, Sections 4.43 and 4.44 focus on the issue of *grouping* observations for a linear regression problem: if data is gathered modulo a rounding process [or if a polyprobit model is to be estimated (Marin and Robert, 2007)], data augmentation (Tanner and Wong, 1987; Robert and Casella, 2004) can recover the original values by simulation, rather than resorting to approximations. Mentions are made of point estimators, but there is unfortunately no connection with decision theory and loss functions in the classical sense (DeGroot, 1970; Berger, 1985). A long section (Section 4.7) deals with rank statistics, containing apparently no connection with Bayesian Statistics, while the final section (Section 4.9) on randomized designs also does not cover the special issue of randomization within Bayesian Statistics (Berger and Wolpert, 1988).

The major components of this chapter in terms of Bayesian theory are an introduction to Laplace's approximation, although not so-called (with an interesting side argument in favor of Jeffreys's priors), some comments on orthogonal parameterisation [understood from an information point of view] and the well-known tramcar example.

### 5.1 Laplace's Approximation

When the number of observations $n$ is large, the posterior distribution can be approximated by a Gaussian centered at the *maximum likelihood estimate* with a *range of order* $n^{-1/2}$. There are numerous instances of the use of Laplace's approximation in Bayesian literature (see, e.g., Berger, 1985; MacKay, 2002), but only with specific purposes oriented toward model choice, not as a generic substitute. Jeffreys derives from this approximation an incentive to treat *the prior probability as uniform* since this is *of no practical importance if the number of observations is large.* His argument is made more precise through the normal approximation,

$$L(\theta|x_1, \ldots, x_n) \approx \tilde{L}(\theta|x) \propto \exp\{-n(\theta - \hat{\theta})^{\mathrm{T}} I(\hat{\theta})(\theta - \hat{\theta})/2\},$$

[27] There is another typo when stating that $\log p/(1-p)$ ranges over $(0, \infty)$.

to the likelihood. [Jeffreys notes that *it is of trivial importance whether* $I(\theta)$ *is evaluated for the actual values or for* the MLE $\hat{\theta}$.] Since the normalization factor is

$$(n/2\pi)^{m/2}|I(\theta)|^{1/2},$$

using Jeffreys's prior $\pi(\theta) \propto |I(\theta)|^{1/2}$ means that the posterior distribution is properly normalized and that *the posterior distribution of* $\theta_i - \hat{\theta}_i$ *is nearly the same* (...) *whether it is taken on data* $\hat{\theta}_i$ *or on* $\theta_i$. This sounds more like a pivotal argument in Fisher's fiducial sense than genuine Bayesian reasoning, but it nonetheless brings an additional argument for using Jeffreys's prior, in the sense that the prior provides the *proper* normalizing factor. Actually, this argument is much stronger than it first looks in that it is at the very basis of the construction of matching priors (Welch and Peers, 1963). Indeed, when considering the proper normalizing constant ($\pi(\theta) \propto |I(\theta)|^{1/2}$), the agreement between the frequentist distribution of the maximum likelihood estimator and the posterior distribution of $\theta$ gets closer by an order of 1.

### 5.2 Outside Exponential Families

When considering distributions that are not from exponential families, sufficient statistics of fixed dimension do not exist, and the MLE is much harder to compute. Jeffreys suggests in Section 4.1 using a minimum $\chi^2$ approximation to overcome this difficulty, an approach which is rarely used nowadays.

A particular example is the poly-$t$ (Bauwens, 1984) distribution

$$\pi(\mu|x_1,\dots,x_s) \propto \prod_{r=1}^{s}\left\{1+\frac{(\mu-x_r)^2}{\nu_r s_r^2}\right\}^{-(\nu_r+1)/2}$$

that *happens when several series of observations yield independent estimates [$x_r$] of the same true value [$\mu$]*. The difficulty with this posterior can now be easily solved via a Gibbs sampler that demarginalizes each $t$ density.

Section 4.3 is not directly related to Bayesian Statistics in that it is considering (best) unbiased estimators, even though the Rao–Blackwell theorem is somehow alluded to. The closest connection with Bayesian Statistics could be that, once summary statistics have been chosen for their availability, a corresponding posterior can be constructed conditional on those statistics.[28] The present equivalent of this proposal would then be to use variational methods (Jaakkola and Jordan, 2000) or ABC techniques (Beaumont, Zhang and Balding, 2002).

An interesting insight is given by the notion of *orthogonal parameters* in Section 4.31, to be understood as the choice of a parameterization such that $I(\theta)$ is diagonal. This orthogonalization is central in the construction of reference priors (Kass, 1989; Tibshirani, 1989; Berger and Bernardo, 1992; Berger, Philippe and Robert, 1998) that are identical to Jeffreys's priors. Jeffreys indicates, in particular, that full orthogonalization is impossible for $m = 4$ and more dimensions.

In Section 4.42 the errors-in-variables model is handled rather poorly, presumably because of computational difficulties: when considering $(1 \le r \le n)$

$$y_r = \alpha\xi + \beta + \varepsilon_r, \quad x_r = \xi + \varepsilon'_r,$$

the posterior on $(\alpha, \beta)$ under standard normal errors is

$$\pi(\alpha,\beta|(x_1,y_1),\dots,(x_n,y_x)) \propto \prod_{r=1}^{n}(t_r^2+\alpha^2 s_r^2)^{-1/2} \cdot \exp\left\{-\sum_{r=1}^{n}\frac{(y_r-\alpha x_r-\beta)^2}{2(t_r^2+\alpha^2 s_r^2)}\right\},$$

which induces a normal conditional distribution on $\beta$ and a more complex $t$-like marginal posterior distribution on $\alpha$ that can still be processed by present-day standards.

Section 4.45 also contains an interesting example of a normal $\mathcal{N}(\mu, \sigma^2)$ sample when *there is a known contribution to the standard error*, that is, when $\sigma^2 > \sigma'^2$ with $\sigma'$ known. In that case, using a flat prior on $\log(\sigma^2 - \sigma'^2)$ leads to the posterior

$$\pi(\mu,\sigma|\bar{x},s^2,n) \propto \frac{1}{\sigma^2-\sigma'^2}\frac{1}{\sigma^{n-1}}\exp\left[-\frac{n}{2\sigma^2}\{(\mu-\bar{x})^2+s^2\}\right],$$

which integrates out over $\mu$ to

$$\pi(\sigma|s^2,n) \propto \frac{1}{\sigma^2-\sigma'^2}\frac{1}{\sigma^{n-2}}\exp\left[-\frac{ns^2}{2\sigma^2}\right].$$

The marginal obviously has an infinite mode (or *pole*) at $\sigma = \sigma'$, but there can be a second (and meaningful) mode if $s^2$ is large enough, as illustrated on Figure 2 (*brown curve*). The outcome is indeed different from using the truncated prior $\pi(\mu, \sigma) \propto 1/\sigma$ (*blue curve*), but to conclude that *the inference using this assessment of the prior probability would be that* $\sigma = \sigma'$ is based once again on the false premise that infinite mass posteriors act like Dirac priors, which is not correct: since $\pi(\sigma|s^2, n)$ does not integrate over $\sigma = \sigma'$, the posterior is simply not defined.[29] In that sense, Jeffreys is thus right in rejecting this prior choice as *absurd*.

[28] A side comment on the first-order symmetry between *the probability of a set of statistics given the parameters and that of the parameters given the statistics* seems to precede the first-order symmetry of the (posterior and frequentist) confidence intervals established in Welch and Peers (1963).

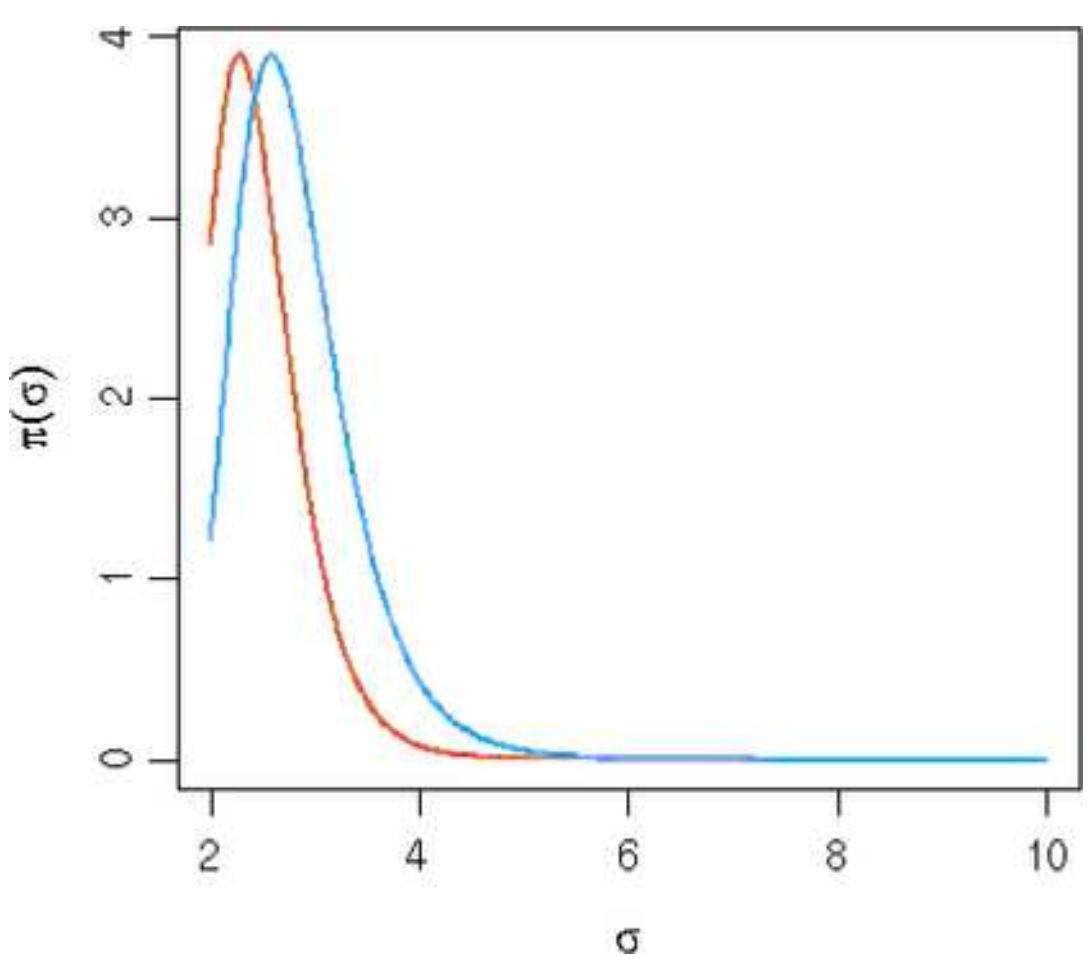

FIG. 2. *Posterior distribution* $\pi(\sigma|s^2, n)$ *for* $\sigma' = \sqrt{2}$, $n = 15$ *and* $ns^2 = 100$, *when using the prior* $\pi(\mu, \sigma) \propto 1/\sigma$ *(blue curve) and the prior* $\pi(\mu, \sigma) \propto 1/\sigma^2 - \sigma'^2$ *(brown curve).*

### 5.3 The Tramcar Problem

This chapter contains (in Section 4.8) the now classical "tramway problem" of Newman, about *a man traveling in a foreign country [who] has to change trains at a junction, and goes into the town, the existence of which he has only just heard. He has no idea of its size. The first thing that he sees is a tramcar numbered* 100. *What can he infer about the number of tramcars in the town? It may be assumed that they are numbered consecutively from* 1 *upwards.*

This is another illustration of the standard noninformative prior for a scale, that is, $\pi(n) \propto 1/n$, where $n$ is the number of tramcars; the posterior satisfies $\pi(n|m = 100) \propto 1/n^2 \mathbb{I}(n \geq 100)$ and

$$\mathbb{P}(n > n_0|m) = \sum_{r=n_0+1}^{\infty} r^{-2} \Big/ \sum_{r=m}^{\infty} r^{-2} \approx \frac{m}{n_0}.$$

Therefore, the posterior median (the justification of which as a Bayes estimator is not included) is approximately $2m$. Although this point is not discussed by Jeffreys, this example is often mentioned in support of the Bayesian approach against the MLE, since the corresponding maximum estimator of $n$ is $m$, always below the true value of $n$, while the Bayes estimator takes a more reasonable value.

## 6. CHAPTER V: SIGNIFICANCE TESTS: ONE NEW PARAMETER

*The essential feature is that we express ignorance of whether the new parameter is needed by taking half the prior probability for it as concentrated in the value indicated by the null hypothesis and distributing the other half over the range possible.*
H. JEFFREYS, *Theory of Probability*, Section 5.0.

This chapter (as well as the following one) is concerned with the central issue of testing hypotheses, the title expressing a focus on the specific case of point null hypotheses: *Is the new parameter supported by the observations, or is any variation expressible by it better interpreted as random?*[30] The construction of Bayes factors as natural tools for answering such questions does require more mathematical rigor when dealing with improper priors than what is found in *Theory of Probability*. Even though it can be argued that Jeffreys's solution (using only improper priors on nuisance parameters) is acceptable via a limiting argument (see also Berger, Pericchi and Varshavsky, 1998, for arguments based on group invariance), the specific and delicate feature of using infinite mass measures would deserve more validation than what is found there. The discussion on the choice of priors to use for the parameters of interest is, however, more rewarding since Jeffreys realizes that (point estimation) Jeffreys's priors cannot be used in this setting (because of their improperness) and that an alternative class of (testing) Jeffreys's priors needs to be introduced. It seems to us that this second type of Jeffreys's priors has been overlooked in the subsequent literature, even though the specific case of the Cauchy prior is often pointed out as a reference prior for testing point null hypotheses involving location parameters.

[29] For an example of a constant MAP estimator, see Robert (2001, Example 4.2).

[30] The formulation of the question restricts the test to embedded hypotheses, even though Section 5.7 deals with normality tests.

### 6.1 Model Choice Formalism

Jeffreys starts by analyzing the question,

> *In what circumstances do observations support a change of the form of the law itself?,*

from a model-choice perspective, by assigning prior probabilities to the models $\mathfrak{M}_i$ that are in competition, $\pi(\mathfrak{M}_i)$ $(i = 1, 2, \ldots)$. He further constrains those probabilities to be *terms of a convergent series*.[31] When checking back in Chapter I (Section 1.62), it appears that this condition is due to the constraint that the probabilities can be normalized to 1, which sounds like an unnecessary condition if dealing with improper priors at the same time.[32] The consequence of this constraint is that $\pi(\mathfrak{M}_i)$ must decrease like $2^{-i}$ or $i^{-2}$ and it thus (a) prevents the use of *equal probabilities* advocated before and (b) imposes an ordering of models.

Obviously, the use of the Bayes factor eliminates the impact of this choice of prior probabilities, as it does for the decomposition of an alternative hypothesis $H_1$ into a series of *mutually irrelevant alternative hypotheses*. The fact that $m$ alternatives are tested at once induces a Bonferroni effect, though, that is not (correctly) taken into account at the beginning of Section 5.04 (even if Jeffreys notes that the Bayes factor is then multiplied by $0.7m$). The following discussion borders more on "ranking and selection" than on testing per se, although the use of Bayes factors with correction factor $m$ or $m^2$ is the proposed solution. It is only at the end of Section 5.04 that the Bonferroni effect of repeated testing is properly recognized, if not correctly solved from a Bayesian point of view.

If the hypothesis to be tested is $H_0 : \theta = 0$, against the alternative $H_1$ that is *the aggregate of other possible values [of θ]*, Jeffreys initiates one of the major advances of *Theory of Probability* by rewriting the prior distribution as a mixture of a point mass in $\theta = 0$ and of a generic density $\pi$ on the range of $\theta$,

$$\pi(\theta) = \tfrac{1}{2}\mathbb{I}_0(\theta) + \tfrac{1}{2}\pi(\theta).$$

This is indeed a stepping stone for Bayesian Statistics in that it explicitly recognizes the need to separate the null hypothesis from the alternative hypothesis within the prior, lest the null hypothesis is not properly weighted once it is accepted. The overall principle is illustrated for a normal setting, $x \sim \mathcal{N}(\theta, \sigma^2)$ (with known $\sigma^2$), so that the Bayes factor is

$$K = \frac{\pi(H_0|x)}{\pi(H_1|x)} \Big/ \frac{\pi(H_0)}{\pi(H_1)} = \frac{\exp\{-x^2/2\sigma^2\}}{\int f(\theta)\exp\{-(x-\theta)^2/2\sigma^2\}\,d\theta}.$$

The numerical calibration of the Bayes factor is not directly addressed in the main text, except via a qualitative divergence from the neutral $K = 1$. Appendix B provides a grading of the Bayes factor, as follows:

- *Grade* 0. $K > 1$. *Null hypothesis supported.*
- *Grade* 1. $1 > K > 10^{-1/2}$. *Evidence against* $H_0$*, but not worth more than a bare mention.*
- *Grade* 2. $10^{-1/2} > K > 10^{-1}$. *Evidence against* $H_0$ *substantial.*
- *Grade* 3. $10^{-1} > K > 10^{-3/2}$. *Evidence against* $H_0$ *strong.*
- *Grade* 4. $10^{-3/2} > K > 10^{-2}$. *Evidence against* $H_0$ *very strong.*
- *Grade* 5. $10^{-2} > K >$. *Evidence against* $H_0$ *decisive.*

The comparison with the $\chi^2$ and $t$ statistics in this appendix shows that a given value of $K$ leads to an increasing (in $n$) value of those statistics, in agreement with Lindley's paradox (see Section 6.3 below).

If there are nuisance parameters $\xi$ in the model (Section 5.01), Jeffreys suggests using the same prior on $\xi$ under both alternatives, $\pi_0(\xi)$, resulting in the general Bayes factor

$$K = \int \pi_0(\xi) f(x|\xi, 0)\,d\xi \Big/ \int \pi_0(\xi)\pi_1(\theta|\xi) f(x|\xi,\theta)\,d\xi\,d\theta,$$

where $\pi_1(\theta|\xi)$ is a conditional density. Note that Jeffreys uses a normal model with Laplace's approximation to end up with the approximation

$$K \approx \frac{1}{\pi_1(\hat\theta|\hat\xi)}\sqrt{\frac{ng_{\theta\theta}}{2\pi}}\exp\left\{-\frac{1}{2}ng_{\theta\theta}\hat\theta^2\right\},$$

[31] The perspective of an infinite sequence of models under comparison is not pursued further in this chapter.

[32] In Jeffreys (1931), Jeffreys puts forward a similar argument *that it is impossible to construct a theory of quantitative inference on the hypothesis that all general laws have the same prior probability* (Section 4.3, page 43). See Earman (1992) for a deeper discussion of this point.

where $\hat{\theta}$ and $\hat{\xi}$ are the MLEs of $\theta$ and $\xi$, and where $g_{\theta\theta}$ is the component of the information matrix corresponding to $\theta$ (under the assumption of strong orthogonality between $\theta$ and $\xi$, which means that the MLE of $\xi$ is identical in both situations). The low impact of the choice of $\pi_0$ on the Bayes factor may be interpreted as a licence to use improper priors on the nuisance parameters despite difficulties with this approach (DeGroot, 1973). An interesting feature of this proposal is that the nuisance parameters are processed independently under both alternatives/models but with the *same* prior, with the consequence that it makes *little difference to $K$ whether we have much or little information about $\theta$*.[33] When the nuisance parameters and the parameter of interest are not orthogonal, the MLEs $\hat{\xi}_0$ and $\hat{\xi}_1$ differ and the approximation of the Bayes factor is now

$$K \approx \frac{\pi_0(\hat{\xi}_0)}{\pi_0(\hat{\xi}_1)} \frac{1}{\pi_1(\hat{\theta}|\hat{\xi}_1)} \sqrt{\frac{n g_{\theta\theta}}{2\pi}} \exp\left\{-\frac{1}{2} n g_{\theta\theta} \hat{\theta}^2\right\},$$

which shows that the choice of $\pi_0$ may have an influence too.

### 6.2 Prior Modeling

In Section 5.02 Jeffreys perceives the difficulty in using an improper prior on the parameter of interest $\theta$ as a normalization problem. If one picks $\pi(\theta)$ or $\pi_1(\theta|\xi)$ as a $\sigma$-finite measure, the Bayes factor $K$ is undefined (rather than *always infinite*, as put forward by Jeffreys when normalizing by $\infty$). He thus imposes $\pi(\theta)$ to be *of any form whose integral converges* (to 1, presumably), ending up in the location case[34] suggesting a Cauchy $\mathcal{C}(0, \sigma^2)$ prior as $\pi(\theta)$.

The first example fully processed in this chapter is the innocuous $\mathcal{B}(n, p)$ model with $H_0 : p = p_0$, which leads to the Bayes factor

$$K = \frac{(n+1)!}{x!(n-x)!} p_0^x (1-p_0)^{n-x} \tag{1}$$

under the uniform prior. While $K = 1$ is recognized as a neutral value, no scaling or calibration of $K$ is mentioned at this stage for reaching a decision about $H_0$ when looking at $K$. The only comment worth noting there is that $K$ is not *very decisive* for small values of $n$: *we cannot get decisive results one way or the other from a small sample* (without adopting a decision framework). The next example still sticks to a compact parameter space, since it deals with the $2 \times 2$ contingency table. The null hypothesis $H_0$ is that of independence between both factors, $H_0 : p_{11} p_{22} = p_{12} p_{21}$. The reparameterization in terms of the margins is

| | 1 | 2 |
|---|---|---|
| 1 | $\alpha\beta + \gamma$ | $\alpha(1-\beta) - \gamma$ |
| 2 | $(1-\alpha)\beta - \gamma$ | $(1-\alpha)(1-\beta) + \gamma$ |

but, in order to simplify the constraint

$$-\min\{\alpha\beta, (1-\alpha)(1-\beta)\} \le \gamma \le \min\{\alpha(1-\beta), (1-\alpha)\beta\},$$

Jeffreys then assumes that $\alpha \le \beta \le 1/2$ via *a mere rearrangement of the table.* In this case, $\pi(\gamma|\alpha, \beta) = 1/\alpha$ over $(-\alpha\beta, \alpha(1-\beta))$. Unfortunately, this assumption (of being able to *rearrange*) is not realistic when $\alpha$ and $\beta$ are unknown and, while the author notes that *in ranges where $\alpha$ is not the smallest, it must be replaced in the denominator [of $\pi(\gamma|\alpha, \beta)$] by the smallest*, the subsequent derivation keeps using the constraint $\alpha \le \beta \le 1/2$ and the denominator $\alpha$ in the conditional distribution of $\gamma$, acknowledging later that *an approximation has been made in allowing $\alpha$ to range from 0 to 1 since $\alpha < \beta < 1/2$.* Obviously, the motivation behind this crude approximation is to facilitate the computation of the Bayes factor,[35] as

$$K \approx \frac{(n_{1\cdot}+1)! n_{2\cdot}! n_{\cdot 1}! n_{\cdot 2}!}{n_{11}! n_{22}! n_{12}! n_{21}! (n+1)!} (n+1)$$

if the data is

| | 1 | 2 | |
|---|---|---|---|
| 1 | $n_{11}$ | $n_{12}$ | $n_{1\cdot}$ |
| 2 | $n_{21}$ | $n_{22}$ | $n_{2\cdot}$ |
| | $n_{\cdot 1}$ | $n_{\cdot 2}$ | $n$ |

The computation of the (true) marginal associated with this prior (under $H_1$) is indeed involved and requires either formal or numerical machine-based integration. For instance, massively simulating from the prior is sufficient to provide this approximation. As shown by Figure 3, the difference between the Monte Carlo approximation and Jeffreys's approximation is not spectacular, even though Jeffreys's approximation appears to be always biased toward larger values, that is, toward the null hypothesis, especially for the values of $K$ larger than 1. In some occurrences, the bias is such that it means acceptance versus rejection, depending on which version of $K$ is used.

[33]The requirement *that $\xi' = \xi$ when $\theta = 0$* (where $\xi'$ denotes the nuisance parameter under $H_1$) seems at first meaningless, since each model is processed independently, but it could signify that the parameterization of both models must be the same when $\theta = 0$. Otherwise, assuming that some parameters are the same under both models is a source of contention within the Bayesian literature.

[34]Note that the section seems to consider only location parameters.

[35]Notice the asymmetry in $n_{1\cdot}$ resulting from the approximation.

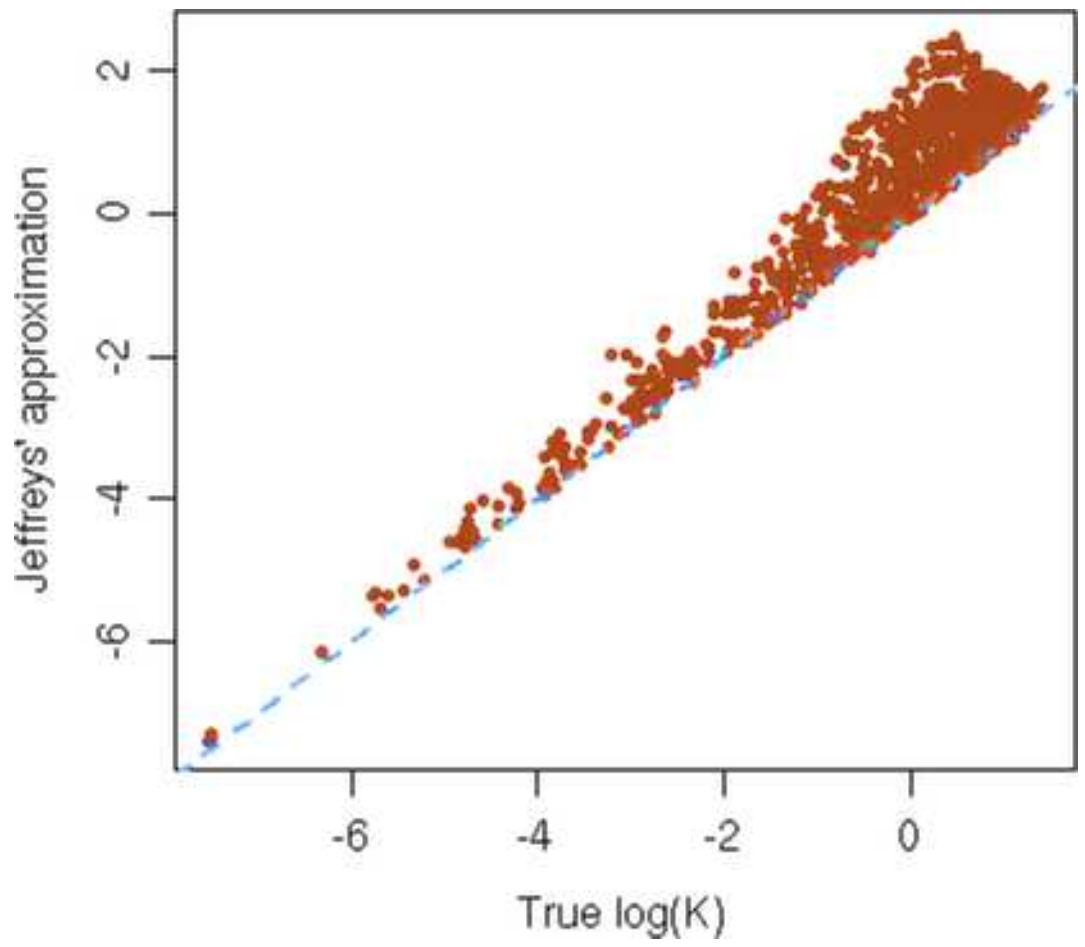

FIG. 3. *Comparison of a Monte Carlo approximation to the Bayes factor for the $2 \times 2$ contingency table with Jeffreys's approximation, based on $10^3$ randomly generated $2 \times 2$ tables and $10^4$ generations from the prior.*

However, if one uses instead a Dirichlet $\mathcal{D}(1,1,1,1)$ prior on the original parameterization $(p_{11},\dots,p_{22})$, the marginal is (up to the multinomial coefficient) the Dirichlet normalizing constant[36]

$$m_1(\mathbf{n}) \propto \frac{D(n_{11}+1,\dots,n_{22}+1)}{D(1,1,1,1)} = 3!\frac{(n+3)!}{n_{11}!n_{22}!n_{12}!n_{21}!},$$

so the (true) Bayes factor in this case is

$$K = \frac{n_{1\cdot}!n_{2\cdot}!n_{\cdot 1}!n_{\cdot 2}!}{((n+1)!)^2}\frac{3!(n+3)!}{n_{11}!n_{22}!n_{12}!n_{21}!} = \frac{n_{1\cdot}!n_{2\cdot}!n_{\cdot 1}!n_{\cdot 2}!}{n_{11}!n_{22}!n_{12}!n_{21}!}\frac{3!(n+3)(n+2)}{(n+1)!},$$

[36]Note that using a Haldane (improper) prior is impossible in this case, since the normalizing constant cannot be eliminated.

which is larger than Jeffreys's approximation. A version much closer to Jeffreys's modeling is based on the parameterization

| | 1 | 2 |
|---|---|---|
| 1 | $\alpha\beta$ | $\gamma(1-\beta)$ |
| 2 | $(1-\alpha)\beta$ | $(1-\gamma)(1-\beta)$ |

in which case $\alpha$, $\beta$ and $\gamma$ are not constrained by one another and a uniform prior on the three parameters can be proposed. After straightforward calculations, the Bayes factor is given by

$$K = (n+1)\frac{n_{\cdot 1}!n_{\cdot 2}!(n_{1\cdot}+1)!(n_{2\cdot}+1)!}{(n+1)!n_{11}!n_{12}!n_{21}!n_{22}!},$$

which is very similar to Jeffreys's approximation since the ratio is $(n_{2\cdot}+1)/(n+1)$. Note that the alternative parameterization based on using

| | 1 | 2 |
|---|---|---|
| 1 | $\alpha\beta$ | $\alpha\gamma$ |
| 2 | $(1-\alpha)(1-\beta)$ | $(1-\alpha)(1-\gamma)$ |

with a uniform prior provides a different answer (with $n_{i\cdot}$'s and $n_{\cdot i}$'s being inverted in $K$). Section 5.12 reprocesses the contingency table with one fixed margin, obtaining very similar outcomes.[37]

In the case of the comparison of two Poisson samples (Section 5.15), $\mathcal{P}(\lambda)$ and $\mathcal{P}(\lambda')$, the null hypothesis is $H_0:\lambda/\lambda' = a/(1-a)$, with $a$ fixed. This suggests the reparameterization

$$\lambda = \alpha\beta, \quad \lambda' = (1-\alpha)\beta',$$

with $H_0:\alpha = a$. This reparameterization appears to be strongly orthogonal in that

$$\begin{aligned} K &= \frac{\int \pi(\beta)a^x(1-a)^{x'}\beta^{x+x'}e^{-\beta}\,d\beta}{\int \pi(\beta)\alpha^x(1-\alpha)^{x'}\beta^{x+x'}e^{-\beta}\,d\beta\,d\alpha} \\ &= \frac{a^x(1-a)^{x'}\int \pi(\beta)\beta^{x+x'}e^{-\beta}\,d\beta}{\int \alpha^x(1-\alpha)^{x'}\,d\alpha \int \pi(\beta)\beta^{x+x'}e^{-\beta}\,d\beta} \\ &= \frac{a^x(1-a)^{x'}}{\int \alpha^x(1-\alpha)^{x'}\,d\alpha} = \frac{(x+x'+1)!}{x!x'!}a^x(1-a)^{x'}, \end{aligned}$$

[37]An interesting example of statistical linguistics is processed in Section 5.14, with the comparison of genders in Welsh, Latin and German, with Freund's *psychoanalytic symbols*, whatever that means!, but the fact that both Latin and German have neuters complicated the analysis so much for Jeffreys that he did without the neuters, apparently unable to deal with $3 \times 2$ tables.

for every prior $\pi(\beta)$, a rather unusual invariance property! Note that, as shown by (1), it also corresponds to the Bayes factor for the distribution of $x$ conditional on $x+x'$ since this is a binomial $\mathcal{B}(x+x',\alpha)$ distribution. The generalization to the Poisson case is therefore marginal since it still focuses on a compact parameter space.

### 6.3 Improper Priors Enter

The bulk of this chapter is dedicated to testing problems connected with the normal distribution. It offers an interesting insight into Jeffreys's processing of improper priors, in that both the infinite mass and the lack of normalizing constant are not clearly signaled as potential problems in the book.

In the original problem of testing the nullity of a normal mean, when $x_1,\dots,x_n \sim \mathcal{N}(\mu,\sigma^2)$, Jeffreys uses a reference prior $\pi_0(\sigma)\propto\sigma^{-1}$ under the null hypothesis and the same reference prior augmented by a proper prior on $\mu$ under the alternative,

$$\pi_1(\mu,\sigma)\propto\frac{1}{\sigma}\pi_{11}(\mu/\sigma)\frac{1}{\sigma},$$

where $\sigma$ is used as a scale for $\mu$. The Bayes factor is then defined as

$$\begin{aligned}K=&\int_0^\infty \sigma^{-n-1}\exp\left\{-\frac{n}{2\sigma^2}(\bar{x}^2+s^2)\right\}d\sigma\\ &\Big/\int_0^\infty\int_{-\infty}^{+\infty}\pi_{11}(\mu/\sigma)\sigma^{-n-2}\\ &\qquad\cdot\exp\left\{-\frac{n}{2\sigma^2}\right.\\ &\qquad\qquad\left.\cdot[(\bar{x}-\mu)^2+s^2]\right\}d\sigma\,d\mu\end{aligned}$$

without any remark on the use of an improper prior in both the numerator and the denominator.[38] There is therefore no discussion about the point of using an improper prior on the nuisance parameters present in both models, that has been defended later in, for example, Berger, Pericchi and Varshavsky (1998) with deeper arguments. The focus is rather on a reference choice for the proper prior $\pi_{11}$. Jeffreys notes that, if $\pi_{11}$ is even, $K=1$ when $n=1$, and he forces the Bayes factor to be zero when $s^2=0$ and $\bar{x}\neq 0$, by a limiting argument that a null empirical variance implies that $\sigma=0$ and thus that $\mu=\bar{x}\neq 0$. This constraint is equivalent to the denominator of $K$ diverging, that is,

$$\int f(v)v^{n-1}\,dv=\infty.$$

A solution[39] that works for all $n\geq 2$ is the Cauchy density, $f(v)=1/\pi(1+v^2)$, advocated as such[40] a reference prior by Jeffreys (while he criticizes the potential use of this distribution for actual data). While the numerator of $K$ is available in closed form,

$$\begin{aligned}&\int_0^\infty \sigma^{-n-1}\exp\left\{-\frac{n}{2\sigma^2}(\bar{x}^2+s^2)\right\}d\sigma\\ &\quad=\left\{\frac{n}{2}(\bar{x}^2+s^2)\right\}^{-n/2}\Gamma(n/2),\end{aligned}$$

this is not the case for the denominator and Jeffreys studies in Section 5.2 some approximations to the Bayes factor, the simplest[41] being

$$K\approx\sqrt{\pi\nu/2}(1+t^2/\nu)^{-(\nu+1)/2},$$

where $\nu=n-1$ and $t=\sqrt{\nu}\bar{x}/s$ (which is the standard $t$ statistic with a constant distribution over $\nu$ under the null hypothesis). Although Jeffreys does not explicitly delve into this direction, this approximation of the Bayes factor is sufficient to expose Lindley's paradox (Lindley, 1957), namely, that the Bayes factor $K$, being equivalent to $\sqrt{\pi\nu/2}\exp\{-t^2/2\}$, goes to $\infty$ with $\nu$ for a fixed value of $t$, thus highlighting the increasing discrepancy between the frequentist and the Bayesian analyses of this testing problem (Berger and Sellke, 1987). As pointed out to us by Lindley (private communication), the paradox is sometimes called the *Lindley–Jeffreys paradox*, because this section clearly indicates that $t$ *increases like* $(\log\nu)^{1/2}$ to keep $K$ constant.

The correct Bayes factor can of course be approximated by a Monte Carlo experiment, using, for instance, samples generated as

$$\sigma^{-2}\sim\mathcal{G}a\left\{\frac{n+1}{2},\frac{ns^2}{2}\right\}\quad\text{and}\quad\mu|\sigma\sim\mathcal{N}(\bar{x},\sigma^2/n).$$

[38] If we extrapolate from earlier remarks by Jeffreys, his justification may be that the *same* normalizing constant (whether or not it is finite) is used in both the numerator and the denominator.

[39] There are obviously many other distributions that also satisfy this constraint. The main drawback of the Cauchy proposal is nonetheless that the scale of 1 is arbitrary, while it clearly has an impact on posterior results.

[40] Cauchy random variables occur in practice as ratios of normal random variables, so they are not completely implausible.

[41] The closest to an explicit formula is obtained just before Section 5.21 as a representation of $K$ through a single integral involving a confluent hypergeometric function.

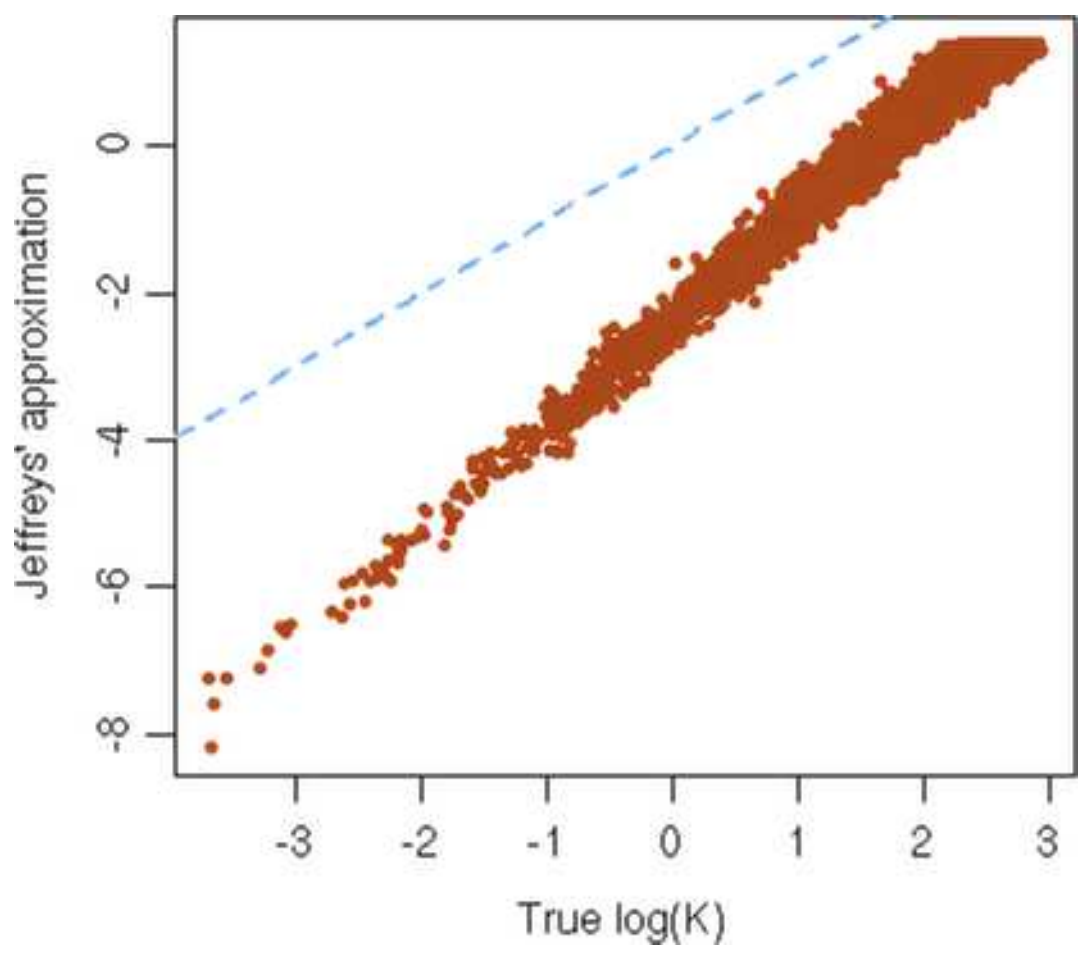

FIG. 4. *Comparison of a Monte Carlo approximation to the Bayes factor for the normal mean problem with Jeffreys's approximation, based on* $5 \times 10^3$ *randomly generated normal sufficient statistics with* $n = 10$ *and* $10^4$ *Monte Carlo simulations of* $(\mu, \sigma)$.

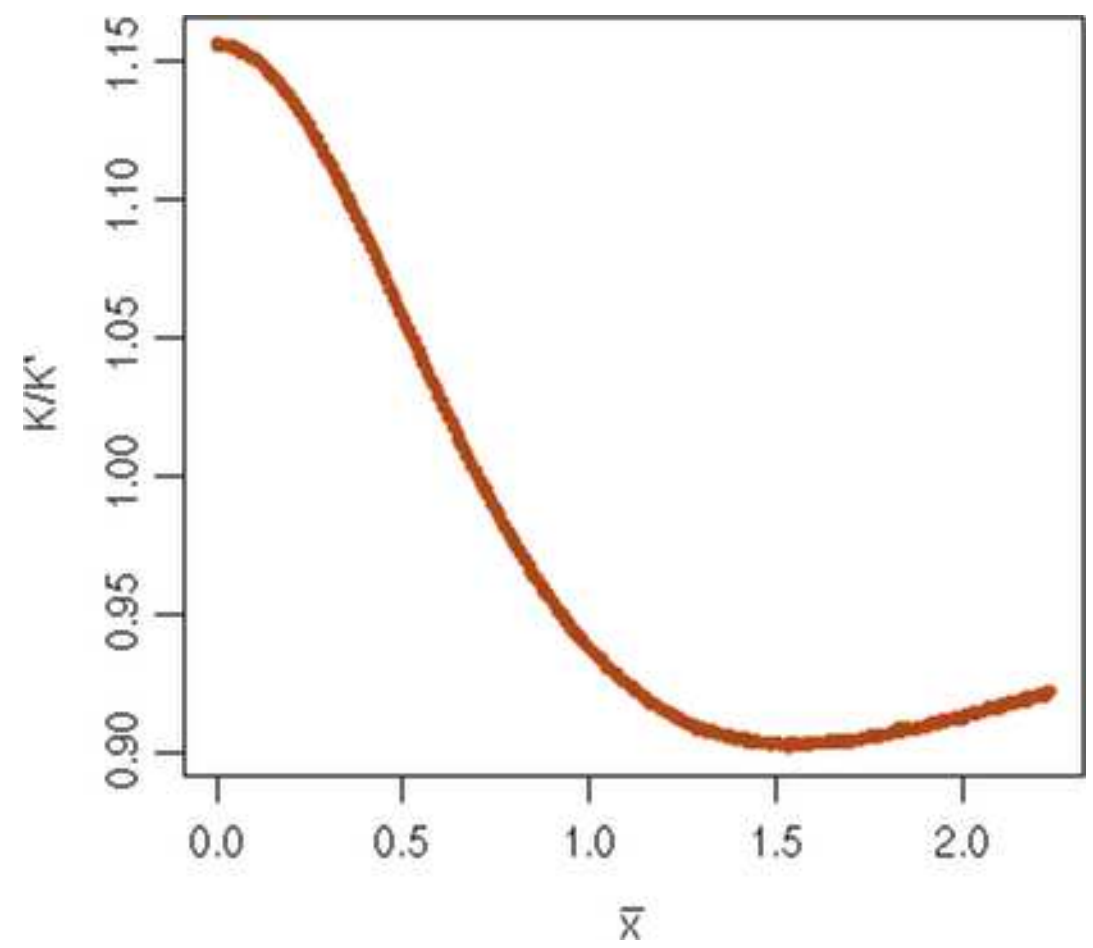

FIG. 5. *Monte Carlo approximation to the Bayes factor for the normal mean problem with known variance, compared with Jeffreys's approximation, based on* $10^6$ *Monte Carlo simulations of* $\mu$*, when* $n = 5$.

The difference between the $t$ approximation and the true value of the Bayes factor can be fairly important, as shown on Figure 4 for $n = 10$. As in Figure 3, the bias is always in the same direction, the approximation penalizing $H_0$ this time. Obviously, as $n$ increases, the discrepancy decreases. (The upper truncation on the cloud is a consequence of Jeffreys's approximation being bounded by $\sqrt{\pi\nu/2}$.)

The Cauchy prior on the mean is also a computational hindrance when $\sigma$ is known: the Bayes factor is then

$$K = \exp\{-n\bar{x}^2/2\sigma^2\} \Big/ \left(\frac{1}{\pi\sigma}\int_{-\infty}^{\infty} \exp\left\{-\frac{n}{2\sigma^2}(\bar{x}-\mu)^2\right\}\frac{d\mu}{1+\mu^2/\sigma^2}\right).$$

In this case, Jeffreys proposes the approximation

$$K \approx \sqrt{2/\pi n}\,\frac{1}{1+\bar{x}^2/\sigma^2},$$

which is then much more accurate, as shown by Figure 5: the maximum ratio between the approximated $K$ and the value obtained by simulation is 1.15 for $n = 5$ and the difference furthermore decreases as $n$ increases.

## 6.4 A Second Type of Jeffreys Priors

In Section 5.3 Jeffreys makes another general proposal for the selection of proper priors under the alternative hypothesis: Noticing that the Kullback divergence is $J(\mu|\sigma) = \mu^2/\sigma^2$ in the normal case above, he deduces that the Cauchy prior he proposed on $\mu$ is equivalent to a flat prior on $\arctan J^{1/2}$:

$$\frac{d\mu}{\pi\sigma(1+\mu^2/\sigma^2)} = \frac{1}{\pi}\frac{dJ^{1/2}}{1+J} = \frac{1}{\pi}d\{\tan^{-1} J^{1/2}(\mu)\},$$

and turns this coincidence into a general rule.[42] In particular, the change of variable from $\mu$ to $J$ is not one-to-one, so there is some technical difficulty linked with this proposal: Jeffreys argues that $J^{1/2}$ should be *taken to have the same sign as* $\mu$ but this is not satisfactory nor applicable in general settings. Obviously, the symmetrization will not always be possible and correcting when the *inverse tangents do not range from* $-\pi/2$ *to* $\pi/2$ can be done in many ways, thus making the idea not fully compatible with the general invariance principle at the core of *Theory of Probability*. Note, however, that Jeffreys's idea of using a functional of the Kullback–Leibler divergence (or of other divergences) as a reference parameterisation for the *new parameter* has many interesting applications. For instance, it is central to the locally conic parameterization used by Dacunha-Castelle and Gassiat (1999) for testing the number of components in mixture models.

In the first case he examines, namely, the case of the contingency table, Jeffreys finds that the corresponding Kullback divergence depends on which margins are fixed (as is well known, the Fisher information matrix is not fully compatible with the Likelihood Principle, see Berger and Wolpert, 1988). Nonetheless, this is an interesting insight that precedes the reference priors of Bernardo (1979): given nuisance parameters, it derives the (conditional) prior on the parameter of interest as the Jeffreys prior for the conditional information. See Bayarri and Garcia-Donato (2007) for a modern extension of this perspective to general testing problems.

[42] We were not aware of this rule prior to reading the book and this second type of Jeffreys's priors, judging from the Bayesian literature, does not seem to have inspired many followers.

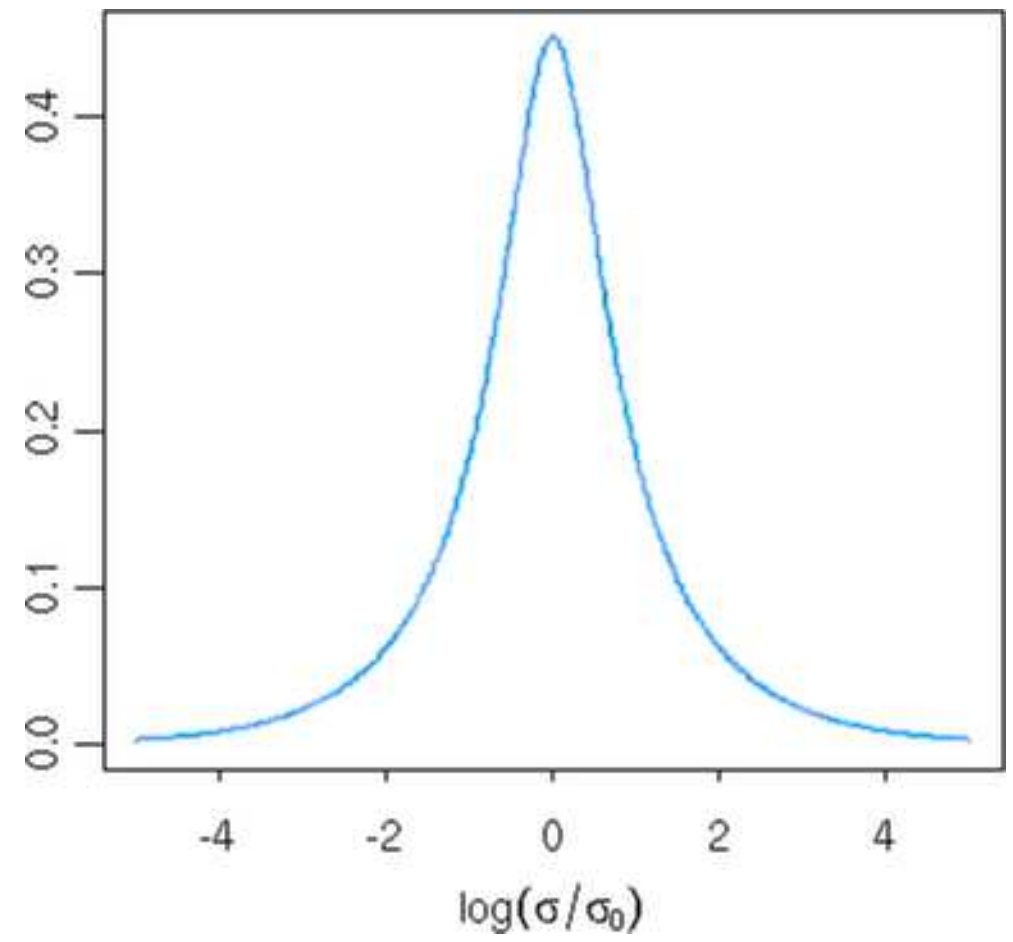

FIG. 6. *Jeffreys's reference density on* $\log(\sigma/\sigma_0)$ *for the test of* $H_0:\sigma=\sigma_0$.

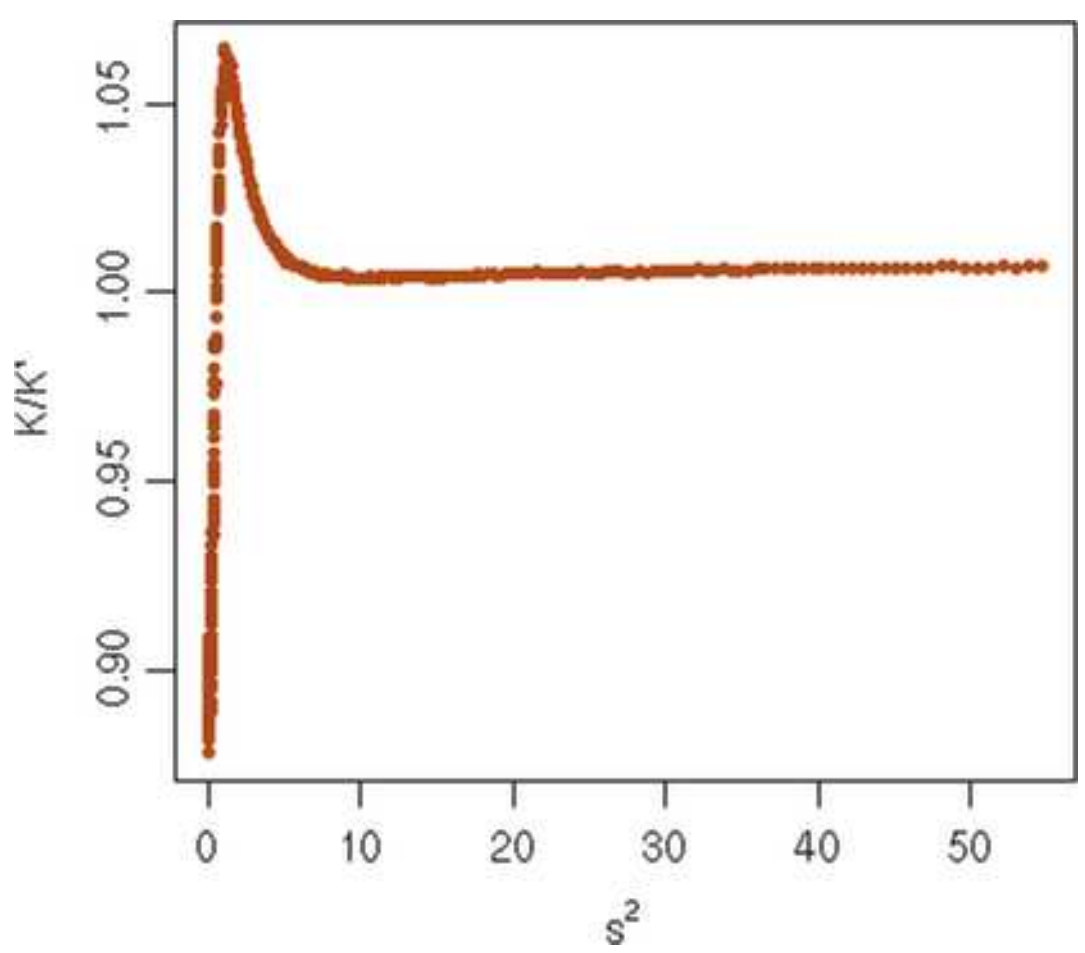

FIG. 7. *Ratio of a Monte Carlo approximation to the Bayes factor for the normal variance problem and of Jeffreys's approximation, when* $n=10$ *(based on* $10^4$ *simulations).*

In the case (Section 5.43) of testing *whether a [normal] standard error has a suggested value* $\sigma_0$ when observing $ns^2 \sim \mathcal{G}a(n/2, \sigma^2/2)$, the parameterization

$$\sigma = \sigma_0 e^{\zeta}$$

leads to (modulo the improper change of variables)

$$J(\zeta) = 2\sinh^2(\zeta) \quad \text{and}$$

$$\frac{1}{\pi}\frac{d\tan^{-1}J^{1/2}(\zeta)}{d\zeta} = \frac{\sqrt{2}\cosh(\zeta)}{\pi\cosh(2\zeta)}$$

as a potential (and overlooked) prior on $\zeta = \log(\sigma/\sigma_0)$.[43] The corresponding Bayes factor is not available in closed form since

$$\int_{-\infty}^{\infty}\frac{\cosh(\zeta)}{\cosh(2\zeta)}e^{-n\zeta}\exp\{-ns^2/2\sigma_0^2 e^{2\zeta}\}\,d\zeta$$
$$= \int_0^{\infty}\frac{1+u^2}{1+u^4}u^n\exp\left\{-\frac{ns^2}{2}u^2\right\}du$$

cannot be analytically integrated, even though a Monte Carlo approximation is readily computed. Figure 7 shows that Jeffreys's approximation,

$$K \approx \sqrt{\pi n/2}\frac{\cosh(2\log s/\sigma_0)}{\cosh(\log s/\sigma_0)}(s/\sigma_0)^n$$
$$\cdot \exp\{n(1-(s/\sigma_0)^2)/2\},$$

is again fairly accurate since the ratio is at worst 0.9 for $n=10$ and the difference decreases as $n$ increases.

The special case of testing a normal correlation coefficient $H_0:\rho=\rho_0$ is not processed (in Section 5.5) via this general approach but, based on arguments connected with (a) the earlier difficulties in the construction of an appropriate noninformative prior (Section 4.7) and (b) the fact that $J$ diverges for the null hypothesis[44] $\rho=\pm 1$, Jeffreys falls back on the uniform $\mathcal{U}(-1,1)$ solution, which is even more convincing in that it leads to an almost closed-form solution

$$K = \frac{2(1-\rho_0^2)^{n/2}/(1-\rho\hat\rho)^{n-1/2}}{\int_{-1}^{1}(1-\rho^2)^{n/2}/(1-\rho\hat\rho)^{n-1/2}\,d\rho}.$$

Note that Jeffreys's approximation,

$$K \approx \left(\frac{2n-1}{\pi}\right)^{1/2}\frac{(1-\rho_0^2)^{n/2}(1-\hat\rho^2)^{(n-3)/2}}{(1-\rho\hat\rho)^{n-1/2}},$$

is quite reasonable in this setting, as shown by Figure 8, and also that the value of $\rho_0$ has no influence on the ratios of the approximations. The extension to two samples in Section 5.51 (for testing whether or not the correlation is the same) is not processed in a symmetric way, with some uncertainty about the validity of the expression for the Bayes factor: a pseudo-common correlation is defined under the alternative in *accordance with the rule that the parameter $\rho$ must appear in the statement of $H_1$*, but normalizing constraints on $\rho$ are not properly assessed.[45]

[43]Note that this is indeed a probability density, whose shape is given in Figure 6, despite the loose change of variables, because a missing 2 cancels with a missing 1/2!

[44]This choice of the null hypothesis is somehow unusual, since, on the one hand, it is more standard to test for no correlation, that is, $\rho=0$, and, on the other hand, having $\rho=\pm 1$ is akin to a unit-root test that, as we know today, requires firmer theoretical background.

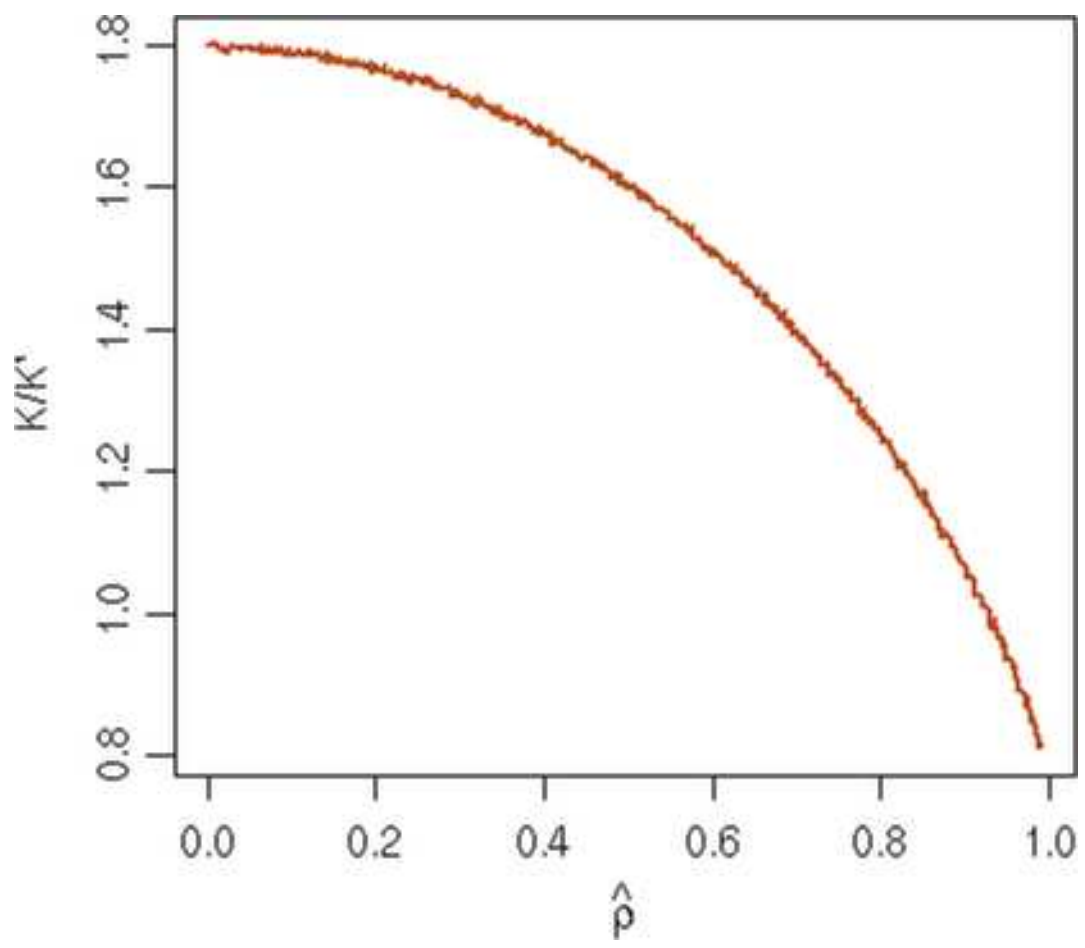

FIG. 8. *Ratio of a Monte Carlo approximation to the Bayes factor for the normal variance problem and of Jeffreys's approximation, when $n = 10$ and $\rho_0 = 0$ (based on $10^4$ simulations).*

A similar approach is adopted for the comparison of two correlation coefficients, with some quasi-hierarchical arguments (see Section 6.5) for the definition of the prior under the alternative. Section 5.6 is devoted to a very specific case of correlation analysis that corresponds to our modern random effect model. A major part of this section argues in favor of the model based on observations in various fields, but the connection with the chapter is the devising of a test for the presence of those random effects. The model is then formalized as normal observations $x_r \sim \mathcal{N}(\mu, \tau^2 + \sigma^2/k_r)$ $(1 \le r \le m)$, where $k_r$ denotes the number of observations within class $r$ and $\tau$ is the variance of the random effect. The null hypothesis is therefore $H_0 : \tau = 0$. Even at this stage, the development is not directly relevant, except for approximation purposes, and the few lines of discussion about the Bayes factor indicate that the (testing) Jeffreys prior on $\tau$ should be in $1/\tau^2$ *for small* $\tau^2$, without further specification. The (numerical) complexity of the problem may explain why Jeffreys differs from his usual processing, although current computational tools obviously allow for a complete processing (modulo the proper choice of a prior on $\tau$) (see, e.g., Ghosh and Meeden, 1984).

[45]To be more specific, a normalizing constant $c$ on the distribution of $\rho_2$ that depends on $\rho$ appears in the closed-form expression of $K$, as, for instance, in equation (14).

Jeffreys also advocates using this principle for testing a normal distribution against alternatives from the Pearson family of distributions in Section 5.7, but no detail is given as to how $J$ is computed and how the Bayes factor is derived. Similarly, for the comparison of the Poisson distribution with the negative binomial distribution in Section 5.8, the form of $J$ is provided for the distance between the two distributions, but the corresponding Bayes factor is only given via a very crude approximation with no mention of the corresponding priors.

In Section 5.9 the extension of the (regular) model to the case of (linear) regression and of variable selection is briefly considered, noticing that (a) for a single regressor (Section 5.91), the problem is exactly equivalent to testing whether or not a normal mean $\mu$ is equal to 0 and (b) for more than one regressor (Section 5.92), the test of nullity of one coefficient can be done conditionally on the others, that is, they can be treated as nuisance parameters under both hypotheses. (The case of linear calibration in Section 5.93 is also processed as a by-product.)

### 6.5 A Foray into Hierarchical Bayes

Section 5.4 explores further tests related to the normal distribution, but Section 5.41 starts with a highly unusual perspective. When testing whether or not the means of two normal samples—with likelihood $L(\mu_1, \mu_2, \sigma)$ proportional to

$$\sigma^{-n_1-n_2} \exp\left\{-\frac{n_1}{2\sigma^2}(\bar{x}_1 - \mu_1)^2 - \frac{n_2}{2\sigma^2}(\bar{x}_2 - \mu_2)^2 - \frac{n_1 s_1^2 + n_2 s_2^2}{2\sigma^2}\right\},$$

—are equal, that is, $H_0 : \mu_1 = \mu_2$, Jeffreys also introduces the value of the common mean, $\mu$, into the alternative. A possible, albeit slightly apocryphal, interpretation is to consider $\mu$ as an hyperparameter that appears both under the null and under the alternative, which is then an incentive to use a single improper prior under both hypotheses (once again because of the lack of relevance of the corresponding pseudo-normalizing constant). But there is still a difficulty with the introduction of three different alternatives with a hyperparameter $\mu$:

$$\mu_1 = \mu \quad \text{and} \quad \mu_2 \neq \mu, \quad \mu_1 \neq \mu \quad \text{and} \quad \mu_2 = \mu,$$
$$\mu_1 \neq \mu \quad \text{and} \quad \mu_2 \neq \mu.$$

Given that $\mu$ has no intrinsic meaning under the alternative, the most logical[46] translation of this multiplication of alternatives is that the three formulations lead to three different priors,

$$\pi_{11}(\mu, \mu_1, \mu_2, \sigma) \propto \frac{1}{\pi} \frac{1}{\sigma^2 + (\mu_2 - \mu)^2} \mathbb{I}_{\mu_1 = \mu},$$
$$\pi_{12}(\mu, \mu_1, \mu_2, \sigma) \propto \frac{1}{\pi} \frac{1}{\sigma^2 + (\mu_1 - \mu)^2} \mathbb{I}_{\mu_2 = \mu},$$
$$\pi_{13}(\mu, \mu_1, \mu_2, \sigma) \propto \frac{1}{\pi^2} \frac{\sigma}{\{\sigma^2 + (\mu_1 - \mu)^2\}\{\sigma^2 + (\mu_2 - \mu)^2\}}.$$

When $\pi_{11}$ and $\pi_{12}$ are written in terms of a Dirac mass, they are clearly identical,

$$\pi_{11}(\mu_1, \mu_2, \sigma) = \pi_{12}(\mu_1, \mu_2, \sigma) \propto \frac{1}{\pi} \frac{1}{\sigma^2 + (\mu_1 - \mu_2)^2}.$$

If we integrate out $\mu$ in $\pi_{13}$, the resulting posterior is

$$\pi_{13}(\mu_1, \mu_2, \sigma) \propto \frac{2}{\pi} \frac{1}{4\sigma^2 + (\mu_1 - \mu_2)^2},$$

whose only difference from $\pi_{11}$ is that the scale in the Cauchy is twice as large. As noticed later by Jeffreys, *there is little to choose between the alternatives*, even though the third modeling makes more sense from a modern, hierarchical point of view: $\mu$ and $\sigma$ denote the location and scale of the problem, no matter which hypothesis holds, with an additional parameter $(\mu_1, \mu_2)$ in the case of the alternative hypothesis. Using a common improper prior under both hypotheses can then be justified via a limiting argument, as in Marin and Robert (2007), because those parameters are common to both models. Seen as such, the Bayes factor

$$\int \sigma^{-n-1} \exp\left\{-\frac{n_1}{2\sigma^2}(\bar{x}_1 - \mu)^2 - \frac{n_2}{2\sigma^2}(\bar{x}_2 - \mu)^2 - \frac{n_1 s_1^2 + n_2 s_2^2}{2\sigma^2}\right\} d\sigma\, d\mu$$
$$\bigg/ \int \frac{\sigma^{-n+1}}{\pi^2} \exp\left\{-\frac{n_1}{2\sigma^2}(\bar{x}_1 - \mu_1)^2 - \frac{n_2}{2\sigma^2}(\bar{x}_2 - \mu_2)^2 - \frac{n_1 s_1^2 + n_2 s_2^2}{2\sigma^2}\right\}$$
$$/(\{\sigma^2 + (\mu_1 - \mu)^2\} \cdot \{\sigma^2 + (\mu_2 - \mu)^2\})\, d\sigma\, d\mu\, d\mu_1\, d\mu_2$$

makes more sense because of the presence of $\sigma$ and $\mu$ on both the numerator and the denominator. While the numerator can be fully integrated into

$$\sqrt{\pi/2n}\Gamma\{(n-1)/2\}(ns_0^2/2)^{-(n-1)/2},$$

where $ns_0^2$ denotes the usual sum of squares, the denominator

$$\int \frac{\sigma^{-n}}{\pi/2} \exp\left\{-\frac{n_1}{2\sigma^2}(\bar{x}_1 - \mu_1)^2 - \frac{n_2}{2\sigma^2}(\bar{x}_2 - \mu_2)^2 - \frac{n_1 s_1^2 + n_2 s_2^2}{2\sigma^2}\right\}$$
$$/(4\sigma^2 + (\mu_1 - \mu_2)^2)\, d\sigma\, d\mu_1\, d\mu_2$$

does require numerical or Monte Carlo integration. It can actually be written as an expectation under the standard noninformative posteriors,

$$\sigma^2 \sim \mathcal{IG}((n-3)/2, (n_1 s_1^2 + n_2 s_2^2)/2),$$
$$\mu_1 \sim \mathcal{N}(\bar{x}_1, \sigma^2/n_1), \quad \mu_2 \sim \mathcal{N}(\bar{x}_2, \sigma^2/n_2),$$

of the quantity

$$\mathfrak{h}(\mu_1, \mu_2, \sigma^2) = \frac{2}{\sqrt{n_1 n_2}} \frac{\Gamma((n-3)/2)\{(n_1 s_1^2 + n_2 s_2^2)/2\}^{-(n-3)/2}}{4\sigma^2 + (\mu_1 - \mu_2)^2}.$$

When simulating a range of values of the sufficient statistics $(n_i, \bar{x}_i, s_i)_{i=1,2}$, the difference between the Bayes factor and Jeffreys's approximation,

$$K \approx 2\left(\frac{\pi}{2} \frac{n_1 n_2}{n_1 + n_2}\right)^{1/2} \cdot \left\{1 + n_1 n_2 n_1 + n_2 \frac{(\bar{x}_1 - \bar{x}_2)^2}{n_1 s_1^2 + n_2 s^2}\right\}^{-(n_1+n_2-1)/2},$$

is spectacular, as shown in Figure 9. The larger discrepancy (when compared to earlier figures) can be attributed in part to the larger number of sufficient statistics involved in this setting.

[46] This does not seem to be Jeffreys's perspective, since he later (in Sections 5.46 and 5.47) adds up the posterior probabilities of those three alternatives, effectively dividing the Bayes factor by 3 or such.

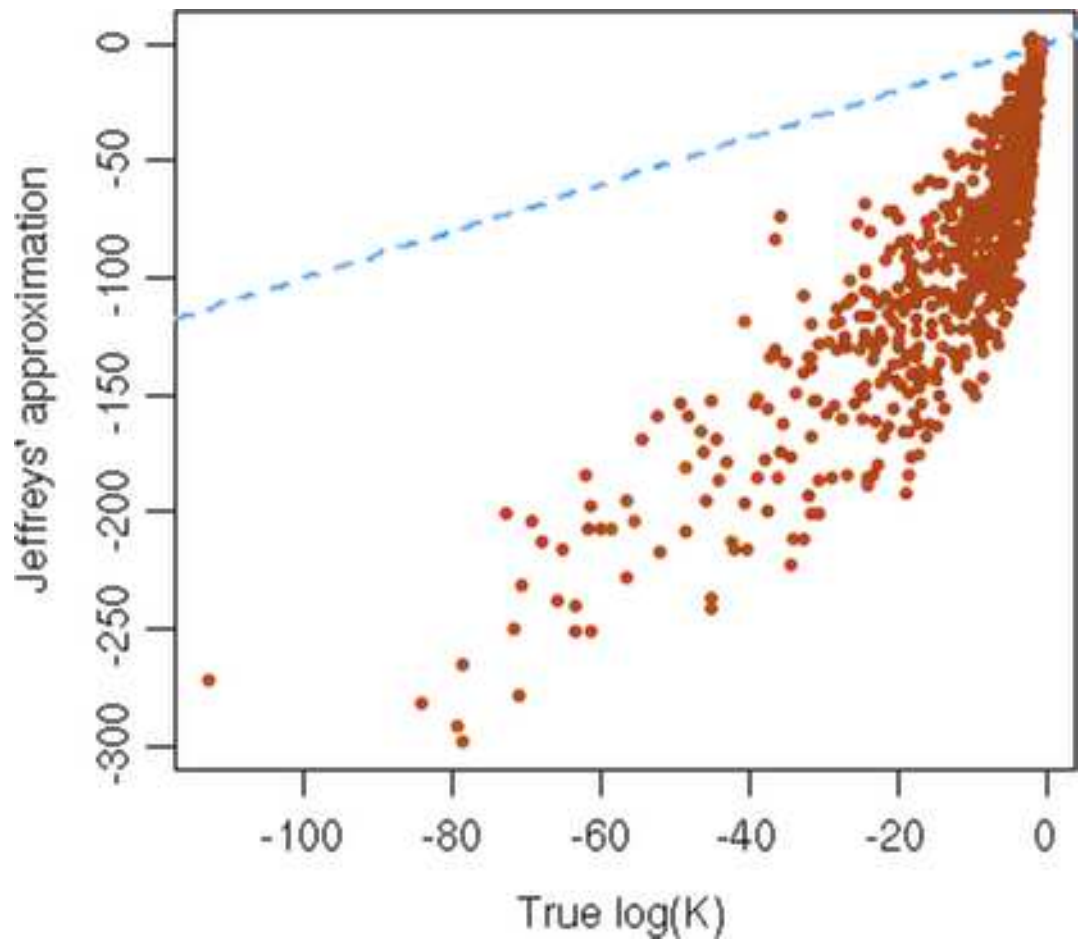

FIG. 9. *Comparison of a Monte Carlo approximation to the Bayes factor for the normal mean comparison problem and of Jeffreys's approximation, corresponding to* $10^3$ *statistics* $(n_i, \bar{x}_i, s_i)_{i=1,2}$ *and* $10^4$ *generations from the noninformative posterior.*

A similar split of the alternative is studied in Section 5.42 when the standard deviations are different under both models, with further simplifications in Jeffreys's approximations to the posteriors (since the $\mu_i$'s are integrated out). It almost seems as if $\bar{x}_1 - \bar{x}_2$ acts as a pseudo-sufficient statistic. If we start from a generic representation with $L(\mu_1, \mu_2, \sigma_1, \sigma_2)$ proportional to

$$\sigma_1^{-n_1}\sigma_2^{-n_2}\exp\left\{-\frac{n_1}{2\sigma_1^2}(\bar{x}_1-\mu_1)^2 - \frac{n_2}{2\sigma_2^2}(\bar{x}_2-\mu_2)^2 - \frac{n_1 s_1^2}{2\sigma_1^2} - \frac{n_2 s_2^2}{2\sigma_2^2}\right\},$$

and if we use again $\pi(\mu, \sigma_1, \sigma_2) \propto 1/\sigma_1\sigma_2$ under the null hypothesis and

$$\pi_{11}(\mu_1, \mu_2, \sigma_1, \sigma_2) \propto \frac{1}{\sigma_1\sigma_2}\frac{1}{\pi}\frac{\sigma_1}{\sigma_1^2 + (\mu_2-\mu_1)^2},$$

$$\pi_{12}(\mu_1, \mu_2, \sigma_1, \sigma_2) \propto \frac{1}{\sigma_1\sigma_2}\frac{1}{\pi}\frac{\sigma_2}{\sigma_2^2 + (\mu_2-\mu_1)^2},$$

$$\pi_{13}(\mu, \mu_1, \mu_2, \sigma_1, \sigma_2) \propto \frac{1}{\sigma_1\sigma_2}\frac{1}{\pi^2}\frac{\sigma_1\sigma_2}{\{\sigma_1^2 + (\mu_1-\mu)^2\}\{\sigma_2^2 + (\mu_2-\mu)^2\}}$$

under the alternative, then, as stated in *Theory of Probability*,

$$\int \sigma_1^{-n_1-1}\sigma_2^{-n_2-1}\exp\left\{-\frac{n_1}{2\sigma_1^2}(\bar{x}_1-\mu)^2 - \frac{n_2}{2\sigma_2^2}(\bar{x}_2-\mu)^2 - \frac{n_1 s_1^2}{2\sigma_1^2} - \frac{n_2 s_2^2}{2\sigma_2^2}\right\}d\mu$$
$$= \sqrt{2\pi/(n_2\sigma_1^2 + n_1\sigma_2^2)}\,\sigma_1^{-n_1}\sigma_2^{-n_2} \cdot \exp\left\{-\frac{(\bar{x}_1-\bar{x}_2)^2}{2(\sigma_1^2/n_1 + \sigma_2^2/n_2)} - \frac{n_1 s_1^2}{2\sigma_1^2} - \frac{n_2 s_2^2}{2\sigma_2^2}\right\},$$

but the computation of

$$\int \exp\left\{-\frac{n_1}{2\sigma_1^2}(\bar{x}_1-\mu_1)^2 - \frac{n_2}{2\sigma_2^2}(\bar{x}_2-\mu)^2\right\} \cdot \frac{2}{2\pi\sigma_2}\frac{d\mu\, d\mu_1}{\sigma_1^2 + (\mu-\mu_1)^2}$$

[and the alternative versions] is not possible in closed form. We note that $\pi_{13}$ corresponds to a distribution on the difference $\mu_1 - \mu_2$ with density equal to

$$\begin{aligned}
&\pi_{13}(\mu_1, \mu_2|\sigma_1, \sigma_2) \\
&\quad = \frac{1}{\pi}((\sigma_1+\sigma_2)(\mu_1-\mu_2)^2 + \sigma_1^3 - \sigma_1^2\sigma_2 - \sigma_1\sigma_2^2 + \sigma_2^3) \\
&\qquad /([(\mu_1-\mu_2)^2 + \sigma_1^2 + \sigma_2^2]^2 - 4\sigma_1^2\sigma_2^2) \\
&\quad = \frac{1}{\pi}\frac{(\sigma_1+\sigma_2)(y^2+\sigma_1^2-2\sigma_1\sigma_2+\sigma_2^2)}{(y^2+(\sigma_1+\sigma_2)^2)(y^2+(\sigma_1-\sigma_2)^2)} \\
&\quad = \frac{1}{\pi}\frac{\sigma_1+\sigma_2}{y^2+(\sigma_1+\sigma_2)^2},
\end{aligned}$$

thus equal to a Cauchy distribution with scale $(\sigma_1 + \sigma_2)$.[47] Jeffreys uses instead a Laplace approximation,

$$\frac{2\sigma_1}{n_1 n_2}\frac{1}{\sigma_1^2 + (\bar{x}_1-\bar{x}_2)^2},$$

to the above integral, with no further justification. Given the differences between the three formulations of the alternative hypothesis, it makes sense to try to compare further those three priors (in our reinterpretation as hierarchical priors). As noted by Jeffreys, *there may be considerable grounds for decision between the alternative hypotheses*. It seems to us (based on the Laplace approximations) that the most sensible prior is the hierarchical one, $\pi_{13}$,

[47] While this result follows from the derivation of the density by integration, a direct proof follows from considering the characteristic function of the Cauchy distribution $\mathcal{C}(0, \sigma)$, equal to $\exp -\sigma|\xi|$ (see Feller, 1971).

in that the scale depends on both variances rather than only one.

An extension of the test on a (normal) standard deviation is considered in Section 5.44 for the *agreement of two estimated standard errors*. Once again, the most straightforward interpretation of Jeffreys's derivation is to see it as a hierarchical modeling, with a reference prior $\pi(\sigma) = 1/\sigma$ on a global scale, $\sigma_1$ say, and the corresponding (testing) Jeffreys prior on the ratio $\sigma_1/\sigma_2 = \exp\zeta$. The Bayes factor (in favor of the null hypothesis) is then given by

$$K = \frac{\sqrt{2}}{\pi} \Big/ \int_{-\infty}^{\infty} \frac{\cosh(\zeta)}{\cosh(2\zeta)} e^{-n_1\zeta} \left(\frac{n_1 e^{2(z-\zeta)} + n_2}{n_2 e^{2z} + n_2}\right)^{-n/2} d\zeta,$$

if $z$ denotes $\log s_1/s_2 = \log \hat\sigma_1/\hat\sigma_2$.

### 6.6 *P*-what?!

Section 5.6 embarks upon a historically interesting discussion on the warnings given by too good a $p$-value: if, for instance, a $\chi^2$ test leads to a value of the $\chi^2$ statistic that is very small, this means (almost certain) incompatibility with the $\chi^2$ assumption just as well as too large a value. (Jeffreys recalls the example of the data set of Mendel that was modified by hand to agree with the Mendelian law of inheritance, leading to *too small a* $\chi^2$ value.) This can be seen as an indirect criticism of the standard tests (see also Section 8 below).

## 7. CHAPTER VI: SIGNIFICANCE TESTS: VARIOUS COMPLICATIONS

*The best way of testing differences from a systematic rule is always to arrange our work so as to ask and answer one question at a time.*
H. Jeffreys, *Theory of Probability*, Section 6.1.

This chapter appears as a *marginalia* of the previous one in that it contains no major advance but rather a sequence of remarks, such as, for instance, an entry on time-series models (see Section 7.2 below). The very first paragraph of this chapter produces a remarkably simple and intuitive justification of the incompatibility between improper priors and significance tests: the mere fact that we are seriously considering the possibility that it is zero may be associated with a presumption that if it is not zero it is probably small.

Then, Section 6.0 discusses the difficulty of settling for an informative prior distribution that takes into account the *actual state of knowledge*. By subdividing the sample into groups, different conclusions can obviously be reached, but this contradicts the Likelihood Principle that the whole data set must be used simultaneously. Of course, this could also be interpreted as a precursor attempt at defining pseudo-Bayes factors (Berger and Pericchi, 1996). Otherwise, as correctly pointed out by Jeffreys, *the prior probability when each subsample is considered is not the original prior probability but the posterior probability left by the previous one*, which is the basic implementation of the Bayesian learning principle. However, even with this correction, the final outcome of a sequential approach is not the proper Bayesian solution, unless posteriors are also used within the integrals of the Bayes factor.

Section 6.5 also recapitulates both Chapters V and VI with general comments. It reiterates the warning, already made earlier, that the Bayes factors obtained via this noninformative approach are usually rarely immensely in favor of $H_0$. This somehow contradicts later studies, like those of Berger and Sellke (1987) and Berger, Boukai and Wang (1997), that the Bayes factor is generally less prone to reject the null hypothesis. Jeffreys argues that, when an alternative is *actually used (...), the probability that it is false is always of order* $n^{-1/2}$, without further justification. Note that this last section also includes the seeds of model averaging: when *a set of alternative hypotheses* (models $\mathfrak{M}_r$) is considered, the predictive should be

$$p(x'|x) = \sum_r p_r(x'|x)\pi(\mathfrak{M}_r|x),$$

rather than conditional on the accepted hypothesis. Obviously, *when $K$ is large, [this] will give almost the same inference as* the selected model/hypothesis.

### 7.1 Multiple Parameters

Although it should proceed from first principles, the extension of Jeffreys's (second) rule for selection priors (see Section 6.4) to several parameters is discussed in Sections 6.1 and 6.2 in a spirit similar to the reference priors of Berger and Bernardo (1992), by pointing out that, if two parameters $\alpha$ and $\beta$ are introduced sequentially against the null hypothesis $H_0 : \alpha = \beta = 0$, testing first that $\alpha \neq 0$ then $\beta \neq 0$ conditional on $\alpha$ does not lead to the same joint prior as the symmetric steps of testing first $\beta \neq 0$ then $\alpha \neq 0$ conditional on $\beta$. In fact,

$$d \arctan J_\alpha^{1/2} d \arctan J_{\beta|\alpha}^{1/2} \neq d \arctan J_\beta^{1/2} d \arctan J_{\alpha|\beta}^{1/2}.$$

Jeffreys then suggests using instead the marginalized version

$$\pi(\alpha, \beta) = \frac{1}{\pi^2} \frac{\frac{dJ_\alpha^{1/2}}{d\alpha}}{1 + J_\alpha} \frac{\frac{dJ_\beta^{1/2}}{d\beta}}{1 + J_\beta},$$

although he acknowledges that there are cases where the symmetry does not make sense (as, for instance, when parameters are not defined under the null, as, e.g., in a mixture setting). He then resorts to Ockham's razor (Section 6.12) to rank those unidimensional tests by stating that *there is a best order of procedures*, although there are cases where such an ordering is arbitrary or not even possible. Section 6.2 considers a two-dimensional parameter $(\lambda, \mu)$ and, switching to polar coordinates, uses a (half-)Cauchy prior on the radius $\rho = \sqrt{\lambda^2 + \mu^2}$ (and a uniform prior on the angle). The Bayes factor for testing the nullity of the parameter $(\lambda, \mu)$ is then

$$\begin{aligned} K = &\int \sigma^{-2n-1} \exp\left\{-\frac{2ns^2 + n(\bar{x}^2 + \bar{y}^2)}{2\sigma^2}\right\} d\sigma \\ &\Big/ \int\!\!\int \frac{1}{\pi^2 \sigma^{2n}} \cdot \exp\{-(2ns^2 + n([\bar{x} - \lambda]^2 + [\bar{y} - \mu]^2))/2\sigma^2\} \frac{d\lambda\, d\mu\, d\sigma}{\rho(\sigma^2 + \rho^2)} \\ = &\, 2^n (n-1)! \{2ns^2 + n(\bar{x}^2 + \bar{y}^2)\}^{-n} \\ &\Big/ \int\!\!\int \frac{1}{\pi^2 \sigma^{2n}} \cdot \exp\left\{-\frac{n}{2\sigma^2} \cdot [2s^2 + \hat{\rho}^2 - 2\rho\hat{\rho}\cos\phi + \rho^2]\right\} \frac{d\phi\, d\rho\, d\sigma}{\rho(\sigma^2 + \rho^2)}, \end{aligned}$$

where $\hat{\rho}^2 = \bar{x}^2 + \bar{y}^2$ and which can only be integrated up to

$$\frac{1}{K} = \frac{2}{\pi} \int_0^\infty \exp\left(-\frac{ns^2 v^2}{2s^2 + \hat{\rho}^2}\right) \cdot {}_1F_1\left\{1 - n, 1, -\frac{n\hat{\rho}^2 v^2}{2(2s^2 + \hat{\rho}^2)}\right\} \frac{dv}{1 + v^2},$$

${}_1F_1$ denoting a confluent hypergeometric function. A similar analysis is conducted in Section 6.21 for a linear regression model associated with a pair of harmonics ($x_t = \alpha \cos t + \beta \sin t + \varepsilon_t$), the only difference being the inclusion of the covariate scales $A$ and $B$ within the prior,

$$\pi(\alpha, \beta|\sigma) = \frac{\sqrt{A^2 + B^2}}{\pi^2 \sqrt{2}} \cdot \frac{\sigma}{\sqrt{\alpha^2 + \beta^2}\{\sigma^2 + (A^2 + B^2)(\alpha^2 + \beta^2)/2\}}.$$

## 7.2 Markovian Models

While the title of Section 6.3 (*Partial and serial correlation*) is slightly misleading, this section deals with an $AR(1)$ model,

$$x_{t+1} = \rho x_t + \tau \varepsilon_t.$$

It is not conclusive with respect to the selection of the prior on $\rho$ given that Jeffreys does not consider the null value $\rho = 0$ but rather $\rho = \pm 1$ which leads to difficulties, if only because there is no stationary distribution in that case. Since the Kullback divergence is given by

$$J(\rho, \rho') = \frac{1 + \rho\rho'}{(1 - \rho^2)(1 - \rho'^2)} (\rho' - \rho)^2,$$

Jeffreys's (testing) prior (against $H_0 : \rho = 0$) should be

$$\frac{1}{\pi} \frac{J^{1/2}(\rho, 0)'}{1 + J(\rho, 0)} = \frac{1}{\pi} \frac{1}{\sqrt{1 - \rho^2}},$$

which is also Jeffreys's regular (estimation) prior in that case.

The (other) correlation problem of Section 6.4 also deals with a Markov structure, namely, that

$$P(x_{t+1} = s | x_t = r) = \begin{cases} \alpha + (1 - \alpha) p_r, & \text{if } s = r, \\ (1 - \alpha) p_s, & \text{otherwise,} \end{cases}$$

the null (independence) hypothesis corresponding to $H_0 : \alpha = 0$. Note that this parameterization of the Markov model means that the $p_r$'s *are* the stationary probabilities. The Kullback divergence being particularly intractable,

$$J = \alpha \sum_{r=1}^m p_r \log\left\{1 + \frac{\alpha}{p_r(1 - \alpha)}\right\},$$

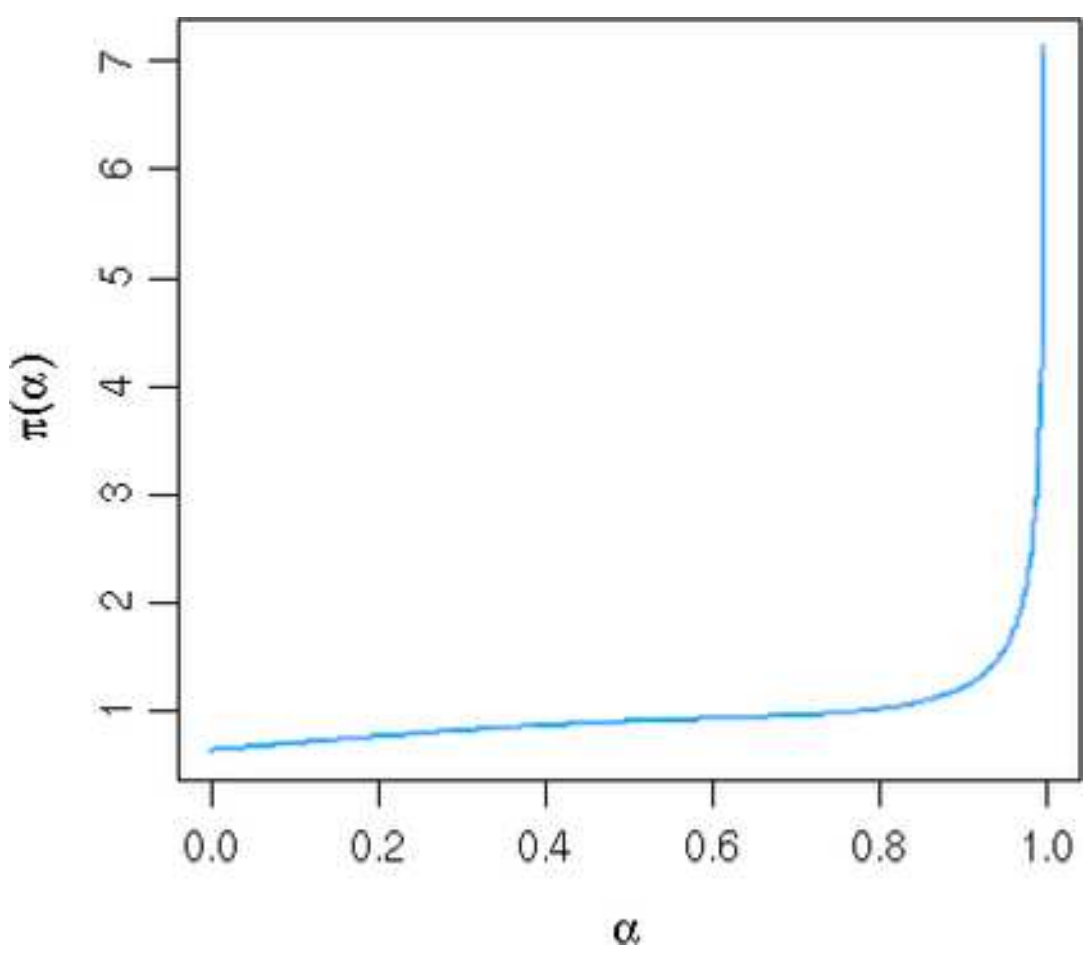

FIG. 10. *Jeffreys's prior of the coefficient $\alpha$ for the Markov model of Section 6.4.*

Jeffreys first produces the approximation

$$J \approx \frac{(m-1)\alpha^2}{1-\alpha}$$

that would lead to the (testing) prior

$$\frac{2}{\pi}\frac{1-\alpha/2}{\sqrt{1-\alpha}(1-\alpha+\alpha^2)}$$

[since the primitive of the above is $-\arctan(\sqrt{1-\alpha}/\alpha)$], but *the possibility of negative*[48] $\alpha$ leads him to use instead a flat prior on the possible range of $\alpha$'s. Note from Figure 10 that the above prior is quite peaked in $\alpha = 1$.

## 8. CHAPTER VII: FREQUENCY DEFINITIONS AND DIRECT METHODS

*An hypothesis that may be true may be rejected because it has not predicted observable results that have not occurred.*

H. JEFFREYS, *Theory of Probability*, Section 7.2.

This short chapter opposes the classical approaches of the time (Fisher's fiducial and likelihood methodologies, Pearson's and Neyman's $p$-values) to the Bayesian principles developed in the earlier chapters. (The very first part of the chapter is a digression on the "frequentist" theories of probability that is not particularly relevant from a mathematical perspective and that we have already addressed earlier. See, however, Dawid, 2004, for a general synthesis on this point.) The fact that Student's and Fisher's analyses of the $t$ statistic coincide with Jeffreys's is seen as an argument in favor both of the Bayesian approach and of the choice of the reference prior $\pi(\mu,\sigma) \propto 1/\sigma$.

The most famous part of the chapter (Section 7.2) contains the often-quoted sentence above, which applies to the criticism of $p$-values, since a decision to reject the null hypothesis is based on the observed $p$-value being in the upper tail of its distribution under the null, even though *nothing but the observed value is relevant*. Given that the $p$-value is a one-to-one transform of the original test statistics, the criticism is maybe less virulent than it appears: Jeffreys still refers to *twice the standard error as a criterion for possible genuineness and three times the standard error for definite acceptance.* The major criticism that this quantity does not account for the alternative hypothesis (as argued, for instance, in Berger and Wolpert, 1988) does not appear at this stage, but only later in Section 7.22. As perceived in *Theory of Probability*, the problem with Pearson's and Fisher's approaches is therefore rather the use of a *convenient* bound on the test statistic as two standard deviations (or on the $p$-value as 0.05). There is, however, an interesting remark that the choice of the hypothesis should eventually be aimed at selecting the best inference, even though Jeffreys concludes that *there is no way of stating this sufficiently precisely to be of any use.* Again, expressing this objective in decision-theoretic terms seems the most natural solution today. Interestingly, the following sentence in Section 7.51 could be interpreted, once again in an apocryphal way, as a precursor to decision theory: *There are cases where there is no positive new parameter, but important consequences might follow if it was not zero*, leading to loss functions mixing estimation and testing as in Robert and Casella (1994).

In Section 7.5 we find a similarly interesting reinterpretation of the classical first and second type errors, computing an integrated error based on the 0–1 loss (even though it is not defined this way) as

$$\int_0^{a_c} f_1(x)\,dx + \int_{a_c}^{\infty} f_0(x)\,dx,$$

where $x$ is the test statistic, $f_0$ and $f_1$ are the marginals under the null and under the alternative, respectively, and $a_c$ is the bound for accepting $H_0$. The

[48]Because of the very specific (unidimensional) parameterization of the Markov chain, using a negative $\alpha$ indeed makes sense.

optimal value of $a_c$ is therefore given by $f_0(a_c) = f_1(a_c)$, which amounts to

$$\pi(H_0|x = a_c) = \pi(H_0^c|x = a_c),$$

that is, $K = 1$ if both hypotheses are equally weighted a priori. This is a completely rigorous derivation of the optimal Bayesian decision for testing, even though Jeffreys does not approach it this way, in particular, because the prior probabilities are not necessarily equal (a point discussed earlier in Section 6.0 for instance). It is nonetheless a fairly convincing argument against $p$-values in terms of *smallest number of mistakes*. More prosaically, Jeffreys briefly discusses in this section the disturbing asymmetry of frequentist tests, when both hypotheses are of the same type: *if we must choose between two definitely stated alternatives, we should naturally take the one that gives the larger likelihood, even though each may be within the range of acceptance of the other.*

## 9. CHAPTER VIII: GENERAL QUESTIONS

*A prior probability used to express ignorance is merely the formal statement of that ignorance.*
H. JEFFREYS, *Theory of Probability*, Section 8.1.

This concluding chapter summarizes the main reasons for using the Bayesian perspective:

1. Prior *and* sampling probabilities are representations of degrees of belief rather than frequencies (Section 8.0). Once again, we believe that this debate[49] is settled today, by considering that probability distributions and improper priors are defined according to the rules of measure theory; see, however, Dawid (2004) for another perspective oriented toward calibration.

2. While prior probabilities are *subjective and cannot be uniquely assessed*, *Theory of Probability* sets a general (objective) principle for the derivation of prior distributions (Section 8.1). It is quite interesting to read Jeffreys's defence of this point when taking into account the fact that this book was setting the point of reference for constructing noninformative priors. *Theory of Probability* does little, however, toward the construction of informative priors by integrating existing prior information (except in the sequential case discussed earlier), recognizing nonetheless the natural discrepancy between two probability distributions conditional on two different data sets. More fundamentally, this stresses that *Theory of Probability* focuses on *prior probabilities used to express ignorance* more than anything else.

3. Bayesian statistics naturally allow for model specification and, as such, do not suffer (as much) from the *neglect of an unforeseen alternative* (Section 8.2). This is obviously true only to some extent: if, in the process of comparing models $\mathfrak{M}_i$ based on an experiment, one very likely model is omitted from the list, the consequences may be severe. On the other hand, and in relation to the previous discussion on the $p$-values, the Bayesian approach allows for alternative models and is thus naturally embedding model specification within its paradigm.[50] The fact that it requires an alternative hypothesis to operate a test is an illustration of this feature.

4. Different theories leading to the same posteriors cannot be distinguished since *questions that cannot be decided by means of observations are best left alone* (Section 8.3). The physicists'[51] concept of *rejection of unobservables* is to be understood as the elimination of *parameters in a law that make no contribution to the results of any observation* or as a version of Ockham's principle, *introducing new parameters only when observations showed them to be necessary* (Section 8.4). See Dawid (1984, 2004) for a discussion of this principle he calls *Jeffreys's Law*.

5. The theory of Bayesian statistics as presented in *Theory of Probability* is consistent in that it provides general rules to construct noninformative priors and to conduct tests of hypotheses (Section 8.6). It is in agreement with the Likelihood Principle and with conditioning on sufficient statistics.[52] It also avoids the use of $p$-values for testing hypotheses by requiring *no empirical hypothesis to be true or false*

[49] Jeffreys argues that the *limit definition was not stated till eighty years later* than Bayes, which sounds incorrect when considering that the Law of Large Numbers was produced by Bernoulli in *Ars Conjectandi*.

[50] The point about *being prepared for occasional wrong decisions* could possibly be related to Popper's notion of *falsifiability*: by picking a specific prior, it is always possible to modify inference toward one's goal. Of course, the divergences between Jeffreys's and Popper's approaches to induction make them quite irreconcilable. See Dawid (2004) for a Bayes–de Finetti–Popper synthesis.

[51] Both paragraphs Sections 8.3 and 8.4 seem only concerned with a physicists' debate, particularly about the relevance of quantum theory.

[52] We recall that Fisher information is not fully compatible with the Likelihood Principle (Berger and Wolpert, 1988).

*a priori*. However, special cases and multidimensional settings show that this theory cannot claim to be completely universal.

6. The final paragraph of *Theory of Probability* states that *the present theory does not justify induction*; what it does *is to provide rules for consistency*. This is absolutely coherent with the above: although the book considers many special cases and exceptions, it does provide a general rule for conducting point inference (estimation) and testing of hypotheses by deriving generic rules for the construction of noninformative priors. Many other solutions are available, but the consistency cannot be denied, while a ranking of those solutions is unthinkable. In essence, *Theory of Probability* has thus mostly achieved its goal of presenting a self-contained theory of inference based on a minimum of assumptions and covering the whole field of inferential purposes.

## 10. CONCLUSION

*It is essential to the possibility of induction that we shall be prepared for occasional wrong decisions.*
H. JEFFREYS, *Theory of Probability*, Section 8.2.

Despite a tone that some may consider as overly critical, and therefore unfair to such a pioneer in our field, this perusal of *Theory of Probability* leaves us with the feeling of a considerable achievement toward the formalization of Bayesian theory and the construction of an objective and consistent framework. Besides setting the Bayesian principle in full generality,

$$\textit{Posterior Probability} \propto \textit{Prior Probability} \cdot \textit{Likelihood},$$

including using improper priors indistinctly from proper priors, the book sets a generic theory for selecting reference priors in general inferential settings,

$$\pi(\theta) \propto |I(\theta)|^{1/2},$$

as well as when testing point null hypotheses,

$$\frac{1}{\pi}\frac{dJ^{1/2}}{1+J} = \frac{1}{\pi} d\{\tan^{-1} J^{1/2}(\theta)\},$$

when $J(\theta) = \operatorname{div}\{f(\cdot|\theta_0), f(\cdot|\theta)\}$ is a divergence measure between the sampling distribution under the null and under the alternative. The lack of a decision-theoretic formalism for point estimation notwithstanding, Jeffreys sets up a completely operational technology for hypothesis testing and model choice that is centered on the Bayes factor. Premises of hierarchical Bayesian analysis, reference priors, matching priors and mixture analysis can be found at various places in the book. That it sometimes lacks mathematical rigor and often indulges in debates that may look superficial today is once again a reflection of the idiosyncrasies of the time: even the ultimate revolutions cannot be built on void and they do need the shoulders of earlier giants to step further. We thus absolutely acknowledge the depth and worth of *Theory of Probability* as a foundational text for Bayesian Statistics and hope that the current review may help in its reassessment.

## ACKNOWLEDGMENTS

This paper originates from a reading seminar held at CREST in March 2008. The authors are grateful to the participants for their helpful comments. Professor Dennis Lindley very kindly provided light on several difficult passages and we thank him for his time, his patience and his supportive comments. We are also grateful to Jim Berger and to Steve Addison for helpful suggestions, and to Mike Titterington for a very detailed reading of a preliminary version of this paper. David Aldrich directed us to his most informative website about Harold Jeffreys and *Theory of Probability*. Parts of this paper were written during the first author's visit to the Isaac Newton Institute in Harold Jeffreys's *alma mater*, Cambridge, whose peaceful working environment was deeply appreciated. Comments from the editorial team of *Statistical Science* were also most helpful.